\documentclass[11pt, leqno]{article}

\usepackage{euscript,amstext,amsmath,amsthm,a4,amssymb}

\newcommand{\comment}[1]{}

\newtheorem{thm}{Theorem}[section]

\newtheorem{prop}[thm]{Proposition}
 \newtheorem{cor}[thm]{Corollary}
\newtheorem{lemma}[thm]{Lemma}
\theoremstyle{remark}
\newtheorem{Rem}[thm]{Remark}
\theoremstyle{definition}
\newtheorem{defn}[thm]{Definition}

\title{A Cheeger-M\"uller theorem for\\ symmetric bilinear torsions}
\author{Guangxiang Su\footnote{Chern Institute of Mathematics \& LPMC, Nankai University,
Tianjin 300071, P.R. China. (sugx@mail.nankai.edu.cn)}\ \ and\
Weiping Zhang\footnote{Chern Institute of Mathematics \& LPMC,
Nankai University, Tianjin 300071, P.R. China.
(weiping@nankai.edu.cn)}}
\date{}

\begin{document}

\maketitle
\begin{abstract}
 We generalize a theorem of Bismut-Zhang, which
extends the Cheeger-M\"uller theorem on Ray-Singer torsion and
Reidemeister torsion, to the case where the flat vector  bundle
over a closed manifold carries a nondegenerate symmetric bilinear
form. As a consequence, we prove the Burghelea-Haller conjecture
which gives an analytic interpretation of the Turaev  torsion. It
thus also provides an analytic interpretation of (not merely the
absolute value of) the Alexander polynomial in knot theory.
\end{abstract}

\renewcommand{\theequation}{\thesection.\arabic{equation}}
\setcounter{equation}{0}

\section{Introduction} \label{s0}

Let $F$ be a unitary flat vector bundle  on  a closed Riemannian
manifold $X$. In \cite{RS}, Ray and Singer defined an analytic
torsion associated to $(X,F)$ and proved that it does not depend
on the Riemannian  metric on $X$. Moreover, they conjectured that
this analytic torsion coincides with the classical  Reidemeister
torsion  defined using a triangulation on $X$ (cf. \cite{Mi}).
This conjecture was later proved in the celebrated papers of
Cheeger \cite{C} and M\"uller \cite{Mu1}. M\"uller generalized
this result in \cite{Mu2} to the case where $F$ is a unimodular
flat vector bundle on $X$. In \cite{BZ1}, inspired by the
considerations of Quillen \cite{Q}, Bismut and
 Zhang reformulated the  above Cheeger-M\"uller theorem as an equality between the Reidemeister and
 Ray-Singer metrics defined on the determinant of
  cohomology, and proved an extension of
 it   to the case of general  flat vector bundles over $X$.
 The method used in \cite{BZ1} is different from those of Cheeger
 and M\"uller in that it makes use of a deformation by Morse functions introduced by Witten
 \cite{W}   on the de Rham complex.

 On the other hand,   Turaev generalizes the concept of
 Reidemeister torsion to a complex valued invariant whose
 absolute value provides the original Reidemeister torsion, with the help of the so-called Euler structure (cf.
 \cite{T}, \cite{FT}). It is natural to ask whether there exists
 an analytic interpretation of this Turaev torsion.

 Recently, there appear  two groups of papers dealing with explicitly this
 question. On   one hand, Braverman and Kappeler \cite{BK1, BK2} define what
 they call ``refined analytic torsion'' for flat vector bundles
 over odd dimensional manifolds, and show that it equals to the
 Turaev torsion up to a multiplication by a complex number of absolute
 value one. On the other hand, Burghelea and Haller \cite{BH1, BH2}, following a suggestion of M\"uller,  define a
 generalized
 analytic torsion associated to a nondegenerate  symmetric bilinear form on  a
 flat vector bundle over an arbitrary dimensional manifold and
 make an explicit conjecture between   this generalized analytic
 torsion and the Turaev torsion.

 Both Braverman-Kappeler and Burghelea-Haller deal with the
 analysis of determinants of non-self-adjoint Laplacians.

 In this paper, we will follow the approach of Burghelea and
 Haller, which is closer in spirit to the approach developed by Bismut-Zhang in \cite{BZ1, BZ2}.

 Let $F$ be  a flat complex vector bundle over an oriented  closed
 manifold. Let $\det H^*(M,F)$ be the determinant line of the
 cohomology with coefficient $F$.

 We make the assumption that $F$ admits a smooth fiberwise
 nondegenerate symmetric bilinear form.\footnote{In general, this might not exist.
 However, as indicated by Burghelea and
 Haller \cite{BH2}, we can form a direct sum of copies of $F$ to make such a symmetric bilinear form
 exists at least on the direct sum.}

 Following Farber-Turaev \cite{FT} and Burghelea-Haller \cite{BH1, BH2},  one constructs naturally a
(nondegenerate)  symmetric bilinear
 form on $\det H^*(M,F)$. This resembles closely with the
 construction of the Ray-Singer metric in \cite{BZ1}, where one
 replaces the symmetric bilinear form by a Hermitian metric on
 $F$. The main difference is that while the Ray-Singer metric is a
 real valued function on elements in $\det H^*(M,F)$, the analytically induced symmetric bilinear form
 generally takes complex values on elements in $\det H^*(M,F)$.

 The main purpose of this paper is to generalize the main result in \cite{BZ1}
 to the current situation. That is
 to say, we establish an explicit comparison result between
  the above analytically  induced symmetric bilinear form on $\det H^*(M,F)$ and
  another one, which is of Reidemeister type,  constructing through a combinatorial way.
  We will state this result in Theorem \ref{t3.1}.

  We will prove this result by the same method as in \cite{BZ1}. That
  is, by making use of the Witten deformation \cite{W} of the de Rham
  complex by a Morse function. However, since we are going to deal
    with complex valued torsion which arises from non-self-adjoint
  Laplacians (the non-self-adjoint property comes from the fact
  that we are dealing with symmetric bilinear forms instead of
  Hermitian metrics), we should take care at each step when we
  will proceed the analytical  arguments  in \cite{BZ1}.  In particular, instead of
  generalizing each step in \cite{BZ1} to the non-self-adjoint case, we will make
  full  use of the results
  in \cite{BZ1} and see what else one needs to do in the current case. It is
  remarkable that everything fits at last to give the desired
  result.

The idea of using the Witten deformation to study symmetric
bilinear torsions was   mentioned before in \cite{BH1}. Moreover,
an important anomaly formula for the analytically constructed
symmetric bilinear forms on  $\det
  H^*(M,F)$  has been
proved in \cite{BH2}.

  A direct consequence of our main result is that if $M$ is of
  vanishing Euler characteristic
  and we
  consider the Euler structure (introduced in \cite{T}) on $M$,  then we can
  prove the Burghelea-Haller conjecture \cite[Conjecture 5.1]{BH2} identifying a modified
  version of the above analytic symmetric bilinear form on $\det
  H^*(M,F)$ with  the Turaev torsion, which is also interpreted as a
symmetric bilinear form on $\det
  H^*(M,F)$.

  Since the Alexander polynomial of a knot in $S^3$ can be expressed by
  certain Turaev torsion (cf. \cite{T} and \cite[Section 7.3]{BH1}),
  our result also provides  a purely  analytic
  interpretation of this famous invariant. This generalizes the
  previously  known  result which expresses the norm of the Alexander
  polynomial by the usual Ray-Singer torsion.

  The rest of this paper is organized as follows. In Section 2, we
  recall the basic definitions of various torsions associated with
  nondegenerate symmetric bilinear forms on a flat vector  bundle, we also state
  an anomaly formula for the analytic torsion associated with
  nondegenerate
  symmetric bilinear forms on a flat vector  bundle. In Section 3,
  we state the main result of this paper and provides a proof of
  it based on several intermediary technical results. Sections 4
  to 9 are
  devoted to the proofs of the intermediary results stated in
  Section 3. In the final Section
  10, we apply the main result proved in Section 3 to prove the
  Burghelea-Haller conjecture \cite[Conjecture 5.1]{BH2} on the analytic interpretation of
  the Turaev torsion. Relations with the Braverman-Kappeler's
  refined analytic torsion \cite{BK1}-\cite{BK4} are also discussed.

  Since we will make substantial use of the results in \cite{BZ1},
  we will refer to \cite{BZ1} for related definitions and
  notations directly  when there will be no confusion.

  The main results of this paper have been announced in \cite{SZ}.

  $
  \ $

\noindent{\bf Acknowledgements} The work of the first named author
was partially supported by the Qiushi Foundation. The work of the
second named author was partially supported by the National
Natural Science Foundation of China.

\section{Symmetric bilinear torsions associated to the de Rham and Thom-Smale complexes} \label{s2}
\setcounter{equation}{0}

 In this section, for a nondegenerate
bilinear symmetric form on a complex flat vector bundle over an
oriented closed manifold, we define  two naturally associated
symmetric bilinear forms on the determinant of the cohomology
$H^*(M,F)$  with coefficient $F$. One constructed in a
combinatorial way through the Thom-Smale complex associated to a
Morse function,  and the other one constructed in an analytic way
through the de Rham complex. An anomaly formula essentially due to
Burghelea-Haller \cite{BH2} of the later will also be recalled.

\subsection{\normalsize  Symmetric bilinear torsion of a finite dimensional complex}
\label{s2.1}

Let $(C^*,\partial)$ be a  finite cochain  complex
\begin{align}\label{2.1}
\left(C^*,\partial\right): 0\longrightarrow
C^0\stackrel{\partial_0}{\longrightarrow}
C^1\stackrel{\partial_1}{\longrightarrow}\cdots\stackrel{\partial_{n-1}}{\longrightarrow}C^n\longrightarrow
0,
\end{align}
where each $C ^i$, $0\leq i\leq n$, is a finite dimensional
complex vector space.

Let
\begin{align}\label{2.2}
H^*\left(C^*,\partial\right)=\bigoplus_{i=0}^{n}H^i\left(C^*,\partial\right),\end{align}
   be the cohomology of $(C^*,\partial )$.

   Let
\begin{align}\label{2.3}
\det \left(C^*,\partial\right)=\bigotimes_{i=0}^{n} \left(\det
C^i\right)^{(-1)^i} ,\end{align}
\begin{align}\label{2.4}
\det H^*\left(C^*,\partial\right)=\bigoplus_{i=0}^{n} \left(\det
H^i\left(C^*,\partial\right)\right)^{(-1)^i}\end{align} be the
determinant lines of $(C^*,\partial)$ and $H^*(C^*,\partial)$
respectively.

It is well-known that there is a canonical isomorphism (cf.
\cite{KM} and \cite[Section 1a)]{BGS})
\begin{align}\label{2.5}
\det \left(C^*,\partial\right)\simeq  \det
H^*\left(C^*,\partial\right).\end{align}

Let each $C^i$, $0\leq i\leq n$, admit a nondegenerate symmetric
bilinear form $b_i$, then by (\ref{2.3}) they induce canonically a
symmetric bilinear form $b_{\det (C^*,\partial)}$ on
$\det(C^*,\partial)$, which in turn, via (\ref{2.5}), induces a
symmetric bilinear form $b_{\det H^*(C^*,\partial )}$ on $\det
H^*(C^*,\partial )$.

\begin{defn}\label{t2.1} (cf. \cite{FT}, \cite{BH1} and \cite{BH2})
We call $b_{\det H^*(C^*,\partial )}$
   the
  symmetric bilinear torsion on $\det H^*(C^*,\partial
)$.
\end{defn}

\begin{Rem}\label{t2.2} If $(C^*,\partial)$ is acyclic, that is, $ H^*(C^*,\partial
)=\{ 0\}$, then  $b_{\det H^*(C^*,\partial )}$ is identified  as a
complex number.
\end{Rem}

Let $A_i$, $0\leq i\leq n$, be an    automorphism of $C^i$. Then
it induces a  symmetric bilinear form $b_i'$  on $C^i$ defined by
\begin{align}\label{2.6}
b_i'(x,y)=b_i(A_ix,A_iy).\end{align} Let $b'_{\det
H^*(C^*,\partial )}$ be the associated   symmetric bilinear
torsion on $\det H^*(C^*,\partial)$.

The following anomaly result is obvious.

\begin{prop} \label{t2.3} The following identity holds,
\begin{align}\label{2.7}
{b'_{\det H^*(C^*,\partial )}\over b_{\det H^*(C^*,\partial
)}}=\prod_{i=0}^n
\left(\det\left(A_i\right)^2\right)^{(-1)^i}.\end{align}
\end{prop}

\subsection{\normalsize  Milnor symmetric bilinear torsion of the Thom-Smale complex}
\label{s2.2}

Let $\  M$ be a    closed    smooth manifold $M$, with $\dim M=n$.
For simplicity, we make the assumption that $M$ is oriented (the
non-orientable case can be treated in exactly the same way, with
obvious modifications).

Let $(F,\nabla^F)$ be a complex flat vector bundle over $M$
carrying   the flat connection $\nabla^F$.  We make the assumption
that $F$ carries a nondegenerate symmetric bilinear form $b^F$.

Let $(F^*,\nabla^{F^*})$ be the dual complex flat vector bundle of
$(F,\nabla^{F})$  carrying   the dual  flat connection
$\nabla^{F^*}$. 

Let $f:M\rightarrow {\bf R}$ be a Morse function. Let $g^{TM}$ be
a Riemannian metric on $TM$ such that the corresponding gradient
vector field $-X=-\nabla f\in \Gamma(TM)$ satisfies the Smale
transversality conditions (cf. \cite{Sm}), that is, the unstable
cells (of $-X$) intersect transversally with the stable cells.

Set
\begin{align}\label{2.8}
B=\{ x\in M; X(x)=0\} .
\end{align}

For any $ {x}\in {B}$, let $W^u( {x})$ (resp. $W^s( {x})$) denote
the unstable (resp. stable) cell at $ {x}$, with respect to $-
{X}$. We also choose an orientation $O_{ {x}}^-$ (resp. $O_{
{x}}^+$) on $W^u( {x})$ (resp. $W^s( {x})$).

Let  $ {x}$, $ {y}\in {B}$ satisfy the Morse index relation ${\rm
ind}( {y})={\rm ind}( {x})-1$, then $\Gamma( {x}, {y})=W^u(
{x})\cap W^s( {y})$ consists of a finite number of integral curves
$\gamma$ of $- {X}$. Moreover, for each $\gamma\in\Gamma( {x},
{y})$, by using the orientations chosen above, on can define a
number $n_\gamma( {x},{y})=\pm 1$ as in \cite[(1.28)]{BZ1}.

If $ {x}\in {B}$, let $[W^u( {x})]$ be the complex line generated
by  $W^u( {x})$. Set
\begin{align}\label{2.9}
C_*(W^u, {F}^*)=\bigoplus_{ {x}\in {B}}[W^u( {x})]\otimes
 {F}^*_{ {x}},
\end{align}
\begin{align}\label{2.10}
C_i(W^u, {F}^*)=\bigoplus_{ {x}\in {B},\ {\rm ind}( {x})=i}[W^u(
{x})]\otimes  {F}^*_{ {x}}.
\end{align}
If $ {x}\in {B}$, the flat vector bundle $ {F}^*$ is canonically
trivialized on $W^u( {x})$. In particular, if $ {x}$, $ {y}\in
{B}$ satisfy  ${\rm ind}( {y})={\rm ind}( {x})-1$, and if
$\gamma\in\Gamma( {x}, {y})$, $f^*\in {F}^*_{ {x}}$, let
$\tau_\gamma(f^*)$ be the parallel transport of $f^*\in {F}^*_{
{x}}$ into $ {F}^*_{ {y}}$ along $\gamma$ with respect to the flat
connection $\nabla^{ {F}^*}$.

Clearly, for any $ {x}\in {B}$, there is only a finite number of $
{y}\in {B}$, satisfying together that ${\rm ind}( {y})={\rm ind}(
{x})-1$ and $\Gamma( {x}, {y})\neq \emptyset$.

If $ {x}\in {B}$, $f^*\in {F}^*_{ {x}}$, set
\begin{align}\label{2.11}
\partial (W^u( {x})\otimes
f^*)=\sum_{ {y}\in {B},\ {\rm ind}( {y})={\rm ind}(x)-1}
\sum_{\gamma\in\Gamma( {x}, {y})}n_\gamma( {x}, {y}) W^u(
{y})\otimes \tau_\gamma(f^*).
\end{align}
Then $\partial$ maps $C_i(W^u, {F}^*)$ into $C_{i-1}(W^u, {F}^*)$.
Moreover, one has
\begin{align}\label{2.12}
\partial^2=0.
\end{align}
That is, $(C_*(W^u, {F}^*),\partial)$ forms a chain complex. We
call it the  Thom-Smale complex associated to $( {M}, F, -X)$.

If $ {x}\in {B}$, let $[W^u( {x})]^*$ be the dual line to $W^u(
{x})$. Let $(C^*(W^u, {F}), {\partial})$ be the complex which is
dual to $(C_*(W^u, {F}^*),\partial)$. For $0\leq i\leq n$, one has
\begin{align}\label{2.13}
C^i(W^u, {F})=\bigoplus_{ {x}\in {B},\ {\rm ind}( {x})=i}[W^u(
{x})]^*\otimes
 {F}_{ {x}}.
\end{align}

Let $W^u( {x})^*\in [W^u( {x})]^*$ be such that $\langle W^u(
{x}), W^u( {x})^*\rangle=1$.

We now introduce a  symmetric bilinear form  on each $[W^u(
{x})]^*\otimes  {F}_{ {x}}$ such that for any $f,\ f'\in {F}_{
{x}}$,
\begin{align}\label{2.14}
\left\langle W^u( {x})^*\otimes f, W^u( {x})^*\otimes
f'\right\rangle =\left\langle f,f'\right\rangle_{b^{ {F}_{ {x}}}}.
\end{align}

 For any $0\leq i\leq n$, let $C^i(W^u, {F})$ carry
the symmetric bilinear form obtained from those   defined in
(\ref{2.14}) so that the splitting (\ref{2.13}) is orthogonal with
respect to it. One verifies that this symmetric bilinear form is
nondegenerate on $C^i(W^u, {F})$.

\begin{defn}\label{t2.4} The symmetric bilinear torsion on the determinant line
of the  cohomology
of the  Thom-Smale cochain complex $(C^*(W^u, {F}), {\partial})$, in
the sense of Definition \ref{t2.1}, is called the  Milnor symmetric
bilnear torsion   associated to $( {M}, F, b^F, -X)$, and is denoted
by $ b^{\cal M}_{( {M}, F, b^F, -X)}$.
\end{defn}

From the anomaly formula (\ref{2.7}), one deduces easily the
following result.

\begin{prop}\label{t2.5} If ${b}^{F}_1$ is another nondegenerate symmetric
bilinear form on the flat vector bundle $F$ over $M$. Let $b^{\cal
M}_{( {M}, F, b^F_1, -X)}$ denote the corresponding symmetric
bilinear torsion on $\det{  H}^*(C^*(W^u, {F}), {\partial})$, then
the following anomaly formula holds,
\begin{align}\label{2.15}
b^{\cal M}_{( {M}, F, b^F_1, -X)}=b^{\cal M}_{( {M}, F, b^F,
-X)}\prod_{x\in B} \det
\left(\left(b^{F|_x}\right)^{-1}b^{F|_x}_1\right)^{(-1)^{{\rm
ind}(x)} }
\end{align}
\end{prop}

\subsection{\normalsize  Ray-Singer symmetric bilinear torsion of the de Rham complex}
\label{s2.3}

We continue the discussion of the previous subsection. However, we
do not use the Morse function and make transversality assumptions.

For any $0\leq i\leq n$, denote
\begin{align}\label{2.16}
\Omega^i( {M}, {F})=\Gamma\left(\Lambda^i(T^* {M})\otimes
 {F}\right),\ \ \ \
\Omega^*( {M}, {F})=\bigoplus_{i=0}^n\Omega^i( {M}, {F}).
\end{align}
Let $d^{ {F}}$ denote the natural exterior differential on
$\Omega^*( {M}, {F})$ induced from $\nabla^{ {F}}$ which maps each
$\Omega^i( {M}, {F})$, $0\leq i\leq n$, into $\Omega^{i+1}( {M},
{F})$.

Let $g^F$ be a Hermitian metric on $F$. The Riemannian metric
$g^{TM}$ and $g^F$ determine a natural  inner product (that is, a
pre-Hilbert space structure) on $\Omega^*( {M}, {F})$  (cf.
\cite[(2.2)]{BZ1} and \cite[(2.3)]{BZ2}).

On the other hand  $g^{TM}$ and the symmetric bilinear form $b^F$
determine together a symmetric bilinear form on $\Omega^*( {M},
{F})$ such that if $u=\alpha f$, $v=\beta g\in\Omega^*( {M}, {F})$
such that $\alpha,\ \beta\in \Omega^*( {M} )$, $f,\ g\in
\Gamma(F)$, then
\begin{align}\label{2.17}
 \langle u,v\rangle_b=\int_M (\alpha\wedge * \beta) b^F(f,g),
\end{align}
where $*$ is the Hodge star operator (cf. \cite{Z}).

Consider the   de Rham complex
\begin{multline}\label{2.18}
\left( \Omega^*( {M}, {F}) ,d^{ {F}}\right):0\rightarrow
 \Omega^0( {M}, {F}) \stackrel{d^{ {F}}}{\rightarrow}
 \Omega^1( {M}, {F})\rightarrow \cdots\\
\stackrel{d^{ {F}}}{\rightarrow}
 \Omega^n( {M}, {F} )\rightarrow 0.
\end{multline}

Let $d^{F *}_b: \Omega^*( {M}, {F})\rightarrow \Omega^*( {M},
{F})$ denote the formal adjoint of $d^{ {F}}$ with respect to the
symmetric bilinear form in (\ref{2.17}). That is, for any $u,\
v\in\Omega^*( {M}, {F})$, one has
\begin{align}\label{2.19}
 \left\langle d^Fu,v\right\rangle_b= \left\langle
 u,d^{F*}_bv\right\rangle_b.
\end{align} Set
\begin{align}\label{2.20}
 {D}_b=d^{ {F}}+d^{ {F}*}_b,\ \ \  {D}_b^2=\left(d^{
{F}}+d^{ {F}*}_b\right)^2=d^{ {F}*}_bd^{ {F}}+ d^{ {F}}d^{
{F}*}_b.
\end{align}
Then the Laplacian $ {D}^2_b$ preserves the ${\bf Z}$-grading of
$\Omega^*( {M}, {F})$.

As was pointed out in \cite{BH1} and \cite{BH2}, $D_b^2$ has the
same principal symbol as the usual Hodge Laplacian (constructed
using the inner product on $\Omega^*(M,F)$ induced from
$(g^{TM},g^F)$) studied for example in \cite{BZ1}.

We collect some well-known facts concerning $D_b^2$ as in
\cite[Proposition 4.1]{BH2}, where the reference \cite{S} is
indicated.

\begin{prop}\label{t2.6} The following properties hold for the Laplacian
$D_b^2$:

(i) The spectrum of $D_b^2$ is discrete. For every $\theta>0$ all
but finitely many points of the spectrum are contained in the
angle $\{ z\in{\bf C}|-\theta<{\rm arg}(z)<\theta\}$;

(ii) If $\lambda$ is in the spectrum of $D_b^2$, then the image of
the associated spectral projection is finite dimensional and
contains  smooth forms only. We refer to this image as the
(generalized) $\lambda$-eigen space of $D_b^2$ and denote it  by
$\Omega^*_{\{\lambda\}}(M,F)$. There exists $N_\lambda\in{\bf N}$
such that
\begin{align}\label{2.21}
  \left.\left(D_b^2-\lambda\right)^{N_\lambda}\right|_{\Omega^*_{\{\lambda\}}(M,F)}=0.
\end{align} We have a $D_b^2$-invariant  $\langle\ ,\
\rangle_b$-orthogonal decomposition
\begin{align}\label{2.22}
  \Omega^* (M,F)=\Omega^*_{\{\lambda\}}(M,F)\oplus \Omega^*_{\{\lambda\}}(M,F)^{\perp } .
\end{align} The restriction of $D_b^2-\lambda$ to $\Omega^*_{\{\lambda\}}(M,F)^{\perp }$ is invertible;

(iii) The decomposition (\ref{2.22}) is invariant under $d^F$ and
$d^{F*}_b$;

(iv) For $\lambda\neq \mu$, the eigen spaces
$\Omega^*_{\{\lambda\}}(M,F)$ and $\Omega^*_{\{\mu\}}(M,F)$ are
$\langle\ ,\ \rangle_b$-orthogonal to each other.
\end{prop}

For any $a\geq 0$, set
\begin{align}\label{2.23}
  \Omega^*_{[0,a]} (M,F)=\bigoplus_{0\leq |\lambda|\leq a}\Omega^*_{\{\lambda\}}(M,F).
\end{align} Let $\Omega^*_{[0,a]} (M,F)^\perp$ denote the $\langle\ ,\ \rangle_b$-orthogonal
complement to $\Omega^*_{[0,a]} (M,F)$.

By \cite[(29)]{BH2} and Proposition \ref{t2.6}, one sees that
$(\Omega^*_{[0,a]} (M,F),d^F)$ forms a finite dimensional complex
whose cohomology equals to that of $(\Omega^*  (M,F),d^F)$.
Moreover, the symmetric bilinear form $\langle\ ,\ \rangle_b$
clearly induces a nondegenerate symmetric bilinear form on each
$\Omega^i_{[0,a]} (M,F)$ with $0\leq i\leq n$. By Definition
\ref{t2.1} one then gets a symmetric bilinear torsion $b_{\det
H^*(\Omega^*_{[0,a]} (M,F),d^F)}$ on $\det H^*(\Omega^*_{[0,a]}
(M,F),d^F)=\det H^*(\Omega^*  (M,F),d^F)$.

For any $0\leq i\leq n$, let $D_{b,i}^2$ be the restriction of
$D_{b}^2$ on $\Omega^i(M,F)$. Then it is shown in \cite{BH2} (cf.
\cite[Theorem 13.1]{S}) that for any $a\geq 0$, the following
regularized zeta determinant is well-defined,
\begin{align}\label{2.24}
   {\det}'\left(D^2_{b,(a,+\infty),i}\right)
   =\exp\left(-\left.{\partial \over\partial s}\right|_{s=0}{\rm Tr}\left[
   \left(\left. D^2_{b,i}\right|_{\Omega^*_{[0,a]} (M,F)^\perp}\right)^{-s}\right]
   \right).
\end{align}

\begin{prop}\label{t2.7}   {\rm (\cite[Proposition 4.7]{BH2})}  The
  symmetric bilinear form on $\det H^*(\Omega^*  (M,F),d^F)$
  defined by
\begin{align}\label{2.25}
 b_{\det
H^*(\Omega^*_{[0,a]} (M,F),d^F)}\prod_{i=0}^n \left(
{\det}'\left(D^2_{b,(a,+\infty),i}\right)\right)^{(-1)^ii}
\end{align} does not depend on the choice of $a\geq 0$.
\end{prop}

\begin{defn}\label{t2.8} The symmetric bilinear form
defined by (\ref{2.25}) is called the Ray-Singer symmetric
bilinear torsion on $\det H^*(\Omega^*  (M,F),d^F)$ and is denoted
by $b^{\rm RS}_{(M,F,g^{TM},b^F)}$.
\end{defn}

\subsection{\normalsize  An anomaly formula for the  Ray-Singer symmetric bilinear  torsion } \label{s2.4}

We continue the discussion of the above subsection.

Let $\theta(F,b^F)\in\Omega^1(M)$ be the Kamber-Tondeur form
defined by (cf. \cite[(4)]{BH2})
\begin{align}\label{2.26}
\theta\left(F,b^F\right)={\rm
Tr}\left[\left(b^F\right)^{-1}\nabla^Fb^F\right].
\end{align}
Then $\theta(F,b^F)$ is a closed one form on $M$ whose cohomology
class  depends only on the homotopy class of $b^F$ (cf.
\cite{BH2}).

Let $\nabla^{TM}$ denote the Levi-Civita connection associated to
the Riemannian metric $g^{TM}$ on $TM$. Let
$R^{TM}=(\nabla^{TM})^2$ be the curvature of $\nabla^{TM}$. Let
$e(TM,\nabla^{TM})\in\Omega^n(M )$ be the associated Euler form
defined by (cf. \cite[(3.17)]{BZ1} and \cite[Chapter 3]{Z})
\begin{align}\label{2.27}
e\left(TM,\nabla^{TM}\right)={\rm Pf}\left(R^{TM}\over
2\pi\right).
\end{align}

Let $g'^{TM}$ be another Riemannian metric on $TM$ and
$\nabla'^{TM}$ be the associated Levi-Civita connection. Let
$\widetilde{e}(TM,\nabla^{TM},\nabla'^{TM})$ be the Chern-Simons
class of $n-1$ smooth forms on $M$, which is defined modulo exact
$n-1$ forms, such that
\begin{align}\label{2.28}
d\widetilde{e}\left(TM,\nabla^{TM},\nabla'^{TM}\right)=e\left(TM,\nabla'^{TM}\right)
-e\left(TM,\nabla^{TM}\right)
\end{align}
(cf. \cite[(4.10)]{BZ1}). Of course, if $n$ is odd,
\begin{align}\label{2.29}
\widetilde{e}\left(TM,\nabla^{TM},\nabla'^{TM}\right)=0.
\end{align}

Let $b'^F$ be another nondegenerate symmetric bilinear form on
$F$.

 Let $b^{\rm RS}_{(M,F,g'^{TM},b'^F)}$
 denote the  Ray-Singer symmetric bilinear torsion  associated to
$g'^{TM}$ and $b'^F$. Then the complex number
$${b^{\rm RS}_{(M,F,g'^{TM},b'^F)}\over
b^{\rm RS}_{(M,F,g^{TM},b^F)}}\in{\bf C}^*$$ is well-defined.

We can now state the anomaly formula, of which an equivalent form
has been proved in \cite[Theorem 4.2]{BH2}, for the Ray-Singer
symmetric bilinear torsion  as follows.

\begin{thm}\label{t2.9} If $b^F$, $b'^F$ lie in the same homotopy class of nondegenerate symmetric bilinear
forms on $F$, then the following identity holds,
\begin{multline}\label{2.30}
{b^{\rm RS}_{(M,F,g'^{TM},b'^F)}\over b^{\rm
RS}_{(M,F,g^{TM},b^F)}}=\exp\left(\int_M
\log\left(\det\left(\left(b^F\right)^{-1}b'^F\right)\right) e\left(TM,\nabla^{TM}\right)\right)\\
{\cdot}\,\exp\left(-\int_M
\theta(F,b'^F)\widetilde{e}\left(TM,\nabla^{TM},\nabla'^{TM}\right)\right).
\end{multline}
In particular, if $\dim M=n$ is odd, then
\begin{align}\label{2.31}
{b^{\rm RS}_{(M,F,g'^{TM},b'^F)}\over b^{\rm
RS}_{(M,F,g^{TM},b^F)}}=1 .
\end{align}
\end{thm}

\begin{Rem}\label{t2.10}
Since $b^F$, $b'^F$ lie in the same homotopy class, one sees that
$\log (\det ( (b^F\ )^{-1}b'^F ) )$ is a  well defined univalent
function on $M$.\end{Rem}

\section{Comparison   between the Ray-Singer and Milnor symmetric bilinear torsions} \label{s3}
\setcounter{equation}{0}

 In this section, we prove the main result
of   this paper, which is an explicit comparison result between
the Ray-Singer and Milnor symmetric bilinear torsions   introduced
in the last section.

The form of the result we will state  formally looks very similar
to a theorem of Bismut-Zhang proved in \cite[Theorem 0.2]{BZ1}, if
one replaces the Hermitian metrics there by the symmetric bilinear
forms. This similarity also reflects in the proof of the main
result here, where we will use as in \cite{BZ1} the Witten
deformation \cite{W} of the de Rham complex by Morse functions.
Moreover, we will make use the analytic techniques developed in
\cite{BZ1} and \cite{BZ2}, some of which go back to the paper of
Bismut-Lebeau \cite{BL}.

Still, since we will deal with non-self-adjoint operators, we have
to generalize  many of the techniques in \cite{BZ1} and \cite{BZ2}
to the current situation. We will point out the differences in due
context.

\subsection{\normalsize   A Cheeger-M\"uller theorem for symmetric bilinear torsions}
\label{s3.1}

We assume that we are in the same situation as in Sections
 \ref{s2.2}-\ref{s2.4}.
By a simple argument of Helffer-Sj\"ostrand \cite[Proposition
5.1]{HS} (cf. \cite[Section 7b)]{BZ1}), we  may and we well assume
that $g^{TM}$ there satisfies the following property  without
altering the  Thom-Smale cochain complex $(C^*(W^u, {F}),
{\partial})$,

 (*): For any $x\in B$, there is a system of coordinates
$y=(y^1,\cdots,y^n)$ centered at $x$ such that near $x$,
\begin{align}\label{3.1}
g^{TM}=\sum_{i=1}^n \left|dy^i\right|^2,\ \ \ f(y)=f(x)-{1\over
2}\sum_{i=1}^{{\rm ind}(x)}\left|y^i\right|^2+{1\over
2}\sum_{i={\rm ind}(x)+1}^{n}\left|y^i\right|^2.
\end{align}

By a result of Laudenbach \cite{La}, $\{W^u(x):x\in B\}$ form a CW
decomposition of $M$.

For any $ {x}\in  {B}$, $ {F}$ is canonically trivialized over
each cell $W^u( {x})$.

Let $ {P}_\infty$ be the de Rham map defined by
\begin{align}\label{3.2}
\alpha\in\Omega^*( {M}, {F})
 \rightarrow
 {P}_\infty\alpha=\sum_{ {x}\in
 {B}}W^u( {x})^*\int_{W^u( {x})}\alpha\in
C^*(W^u, {F}).
\end{align}

 By the Stokes theorem,   one has
\begin{align}\label{3.3}
 {\partial} {P}_\infty= {P}_\infty
d^{ {F}}.
\end{align}
Moreover, it is shown in \cite{La} that $P_\infty$ is a  ${\bf
Z}$-graded   quasi-isomorphism, inducing a canonical  isomorphism
\begin{align}\label{3.4}
 {P}^{ H}_\infty:{  H}^*
 \left(\Omega^*( {M}, {F}),d^{ {F}}\right)\rightarrow
{  H}^*\left(C^*\left(W^u, {F}\right), {\partial}\right),
\end{align}
which in turn induces a natural isomorphism between the
determinant lines,
\begin{align}\label{3.5}
 {P}^{\det {  H}}_\infty:\det{  H}^*\left(\Omega^*\left( {M}, {F}\right),d^{ {F}}
 \right)\rightarrow
\det{ H}^*\left(C^*\left(W^u, {F}\right), {\partial}\right).
\end{align}

Now let $h^{TM}$ be an arbitrary smooth metric on $TM$.

By Definition \ref{t2.8}, one has an associated  Ray-Singer
symmetric bilinear torsion $ b^{\rm RS}_{( {M},F,h^{TM},b^F)} $ on
$\det{  H}^* (\Omega^* ( {M}, {F} ),d^{ {F}}
  )$.
From (\ref{3.5}), one gets a well-defined symmetric bilinear form
\begin{align}\label{3.6}
 {P}^{\det {  H}}_\infty \left(
b^{\rm RS}_{( {M},F,h^{TM},b^F)}\right)
\end{align}
on $\det{ H}^* (C^* (W^u, {F} ), {\partial} )$.

On the other hand, by Definition \ref{t2.4}, one has a
well-defined  Milnor symmetric bilinear torsion
 $
b^{\cal M}_{( {M},F,b^F,-X)}$ on $\det{ H}^*(C^*(W^u, {F}),
{\partial})$,
 where $X=\nabla f$ is the gradient vector field of $f$ associated
to $g^{TM}$.

Let $\psi(TM,\nabla^{TM})$ be the Mathai-Quillen current
(\cite{MQ}) over $TM$, associated to $h^{TM}$, defined  in
\cite[Definition 3.6]{BZ1}. As indicated in \cite[Remark
3.8]{BZ1}, the pull-back current $X^*\psi(TM,\nabla^{TM})$ is
well-defined over $M$.

The main result of this paper, which generalizes  \cite[Theorem
0.2]{BZ1} to the case where $F$ admits a nondegenerate symmetric
bilinear form, can be stated as follows.

\begin{thm}\label{t3.1} The following identity in ${\bf C}$ holds,
\begin{align}\label{3.7}
 { {P}^{\det {  H}}_\infty \left(
b^{\rm RS}_{( {M},F,h^{TM},b^F)}\right)\over b^{\cal M}_{(
{M},F,b^F,-X)}   }=\exp\left(-
\int_M\theta\left(F,b^F\right)X^*\psi\left(TM,\nabla^{TM}\right)\right).
\end{align}
\end{thm}

\begin{Rem}\label{t3.2} By proceeding similarly as in
\cite[Section 7b)]{BZ1}, in order to prove (\ref{3.7}), we may
well assume that $h^{TM}=g^{TM}$. Moreover, we may assume that
 $b^F$, as well as the Hermitian metric $g^F$ on $F$, are
  flat on an open neighborhood of the zero set $B$ of
$X$. From now on, we will make these assumptions.
\end{Rem}

\subsection{\normalsize   Some intermediate results}
\label{s3.2}

We assume that the assumptions made in Remark \ref{t3.2} hold.

For any $T\in {\bf R}$, let $b^F_T$ be the deformed symmetric
bilinear form on $F$ defined by
\begin{align}\label{3.8}
  b^F_T(u,v)=e^{-2Tf}b^F(u,v).
\end{align}
Let $d^{F*}_{b_T}$ be the associated formal adjoint in the sense
of (\ref{2.19}). Set
\begin{align}\label{3.9}
 {D}_{b_T}=d^{ {F}}+d^{ {F}*}_{b_T},\ \ \  {D}_{b_T}^2=\left(d^{
{F}}+d^{ {F}*}_{b_T}\right)^2=d^{ {F}*}_{b_T}d^{ {F}}+ d^{ {F}}d^{
{F}*}_{b_T}.
\end{align}

Let $\Omega^*_{[0,1],T}(M,F)$ be defined as in (\ref{2.23}) with
respect to $ {D}_{b_T}^2$, and let $\Omega^*_{[0,1],T}(M,F)^\perp$
be the corresponding $\langle\ ,\ \rangle_{b_T}$-orthogonal
complement.

Let $P_T^{[0,1]}$ be the orthogonal projection from
$\Omega^*(M,F)$ to $\Omega^*_{[0,1],T}(M,F)$ with respect to the
inner product determined by $g^{TM}$ and $g_T^F=e^{-2Tf}g^F$.  Set
$P_T^{(1,+\infty)}={\rm Id}-P_T^{[0,1]}$.

Following \cite[(7.13)-(7.15)]{BZ1}, we introduce the notations
\begin{align}\label{3.10}
\chi(F)=\sum_{i=0}^{\dim M}(-1)^i\dim H^i(M,F)={\rm
rk}(F)\sum_{x\in B}(-1)^{{\rm ind}(x)},
\end{align}
$$ {\chi}'(F)={\rm rk}(F)\sum_{x\in B}(-1)^{{\rm
ind}(x)}{\rm ind}(x)={\rm rk}(F)\sum_{i=0}^n(-1)^iiM_i,$$ $$ {\rm
Tr}^B_s[f]=\sum_{x\in B}(-1)^{{\rm ind}(x)}f(x),$$ where for any
$0\leq i\leq n$, $M_i$ is the number of $x\in B$ of index $i$.

Let $N$ be the number operator on $\Omega^*(M,F)$ acting on
$\Omega^i(M,F)$ by multiplication by $i$.

We now state several intermediate results whose proofs will be
given later in Sections 4 to 9.

\begin{thm}\label{t3.3} {\rm (Compare with \cite[Theorem 7.6]{BZ1})}
 Let $P_T^{[0,1]}$ be the restriction of $P_\infty$ on
$\Omega^*_{[0,1],T}(M,F)$, let $P_T^{[0,1],\det H}$ be the induced
isomorphism on cohomology, then the following identity holds,
\begin{align}\label{3.11}
\lim_{T\rightarrow +\infty} {P_T^{[0,1],\det H}\left(b_{\det
H^*(\Omega^*_{[0,1],T} (M,F),d^F)}\right)\over b^{\cal M}_{(
{M},F,b^F,-X)}  }\left({T\over \pi}\right)^{{n\over 2}\chi(F)-
{\chi}'(F)}\exp\left(2\,{\rm rk}(F){\rm Tr}^B_s[f]T\right)
\end{align}$$ =1.$$
\end{thm}

\begin{thm}\label{t3.4} {\rm (Compare with \cite[Theorem 7.8]{BZ1})} For any $t>0$,
\begin{align}\label{3.12}
\lim_{T\rightarrow +\infty}  {\rm Tr}_s\left[
N\exp\left(-tD_{b_T}^2\right)P_T^{(1,+\infty)}\right]=0.
\end{align}
Moreover, for any  $d>0$ there exist $c>0$, $C>0$ and $T_0\geq 1$
such that for any $t\geq d$ and  $T\geq T_0$,
\begin{align}\label{3.13}
 \left|{\rm Tr}_s\left[N\exp\left(-tD_{b_T}^2\right)P_T^{(1,+\infty)}\right]\right|\leq c\exp(-Ct).
\end{align}
\end{thm}

\begin{thm}\label{t3.5} {\rm (Compare with \cite[Theorem 7.9]{BZ1})}
For $T\geq 0$ large enough, then
\begin{align}\label{3.14}
  \dim \Omega^{i}_{[0,1],T}(M,F)={\rm rk}(F)\, M_i.
\end{align}
Also,
\begin{align}\label{3.15}
 \lim_{T\rightarrow +\infty}  {\rm Tr}\left[D_{b_T}^2P_T^{[0,1]}\right]=0 .
\end{align}
\end{thm}

For the next results, we will make use the same notation for
Clifford multiplications and Berezin integrals as in \cite[Section
4]{BZ1}.

\begin{thm}\label{t3.6} {\rm (Compare with \cite[Theorem 7.10]{BZ1})}
As $t\rightarrow 0$, the following
identity holds,
\begin{align}\label{3.16}
   {\rm Tr}_s\left[N\exp\left(-tD_{b_T}^2\right)\right]={n\over 2}\chi(F)+O(t)\ \ if\ n\ {\rm is\ even},
\end{align}
$$
  ={\rm rk}(F)\int_M\int^BL\exp\left(-{\dot{R}^{TM}\over 2}\right)
  {1\over \sqrt{t}}+O\left(\sqrt{t}\right)\ \ if\ n \ {\rm is\ odd}.
$$
\end{thm}

\begin{thm}\label{t3.7} {\rm (Compare with \cite[Theorem A.1]{BZ2})}
There exist $0<\alpha\leq 1$, $C>0$ such that for any $0<t\leq
\alpha$, $0\leq T\leq {1\over t}$, then
\begin{multline}\label{3.17}
 \left|  {\rm Tr}_s\left[N\exp\left(-\left(tD_{b }+T\widehat{c}
 (\nabla f)\right)^2\right)\right]-
 {1\over t}\int_M\int^BL\exp\left(-B_{T^2}\right){\rm rk}(F)\right.\\
\left.  -{T\over
2}\int_M\theta\left(F,b^F\right)\int^B\widehat{df}\exp\left(-B_{T^2}\right)-
  {n\over 2}\chi(F)\right|\leq Ct.
\end{multline}
\end{thm}

\begin{thm}\label{t3.8} {\rm (Compare with \cite[Theorem A.2]{BZ2})} For any
$T>0$, the following identity holds,
\begin{multline}\label{3.18}
 \lim_{t\rightarrow 0}  {\rm Tr}_s\left[N\exp\left(-\left(tD_{b }+{T\over t}\widehat{c}
 (\nabla f)\right)^2\right)\right]\\ =
 {1\over 1-e^{-2T}}\left(\left(1+e^{-2T}\right)  {\chi}'(F)-
 ne^{-2T}\chi(F)\right).
\end{multline}
\end{thm}

\begin{thm}\label{t3.9} {\rm (Compare with \cite[Theorem A.3]{BZ2})}
There exist $\alpha\in(0,1]$, $c>0$, $C>0$ such that for any $t\in
(0,\alpha]$, $T\geq 1$, then
\begin{align}\label{3.19}
 \left|  {\rm Tr}_s\left[N\exp\left(-\left(tD_{b }+{T\over t}\widehat{c}
 (\nabla f)\right)^2\right)\right]-
  {\chi}'(F)\right|\leq c\exp(-CT)
  .
\end{align}
\end{thm}

Clearly, we may and we will  assume that the number $\alpha>0$ in
Theorems \ref{3.7} and \ref{3.9} have been chosen to be the same.

\subsection{\normalsize   Proof of Theorem \ref{t3.1}}
\label{s3.3}

First of all, by the anomaly formula (\ref{2.30}), for any $T\geq
0$, one has
\begin{multline}\label{3.20}
{P_T^{[0,1],\det H}\left(b_{\det
H^*\left(\Omega^*_{[0,1],T}(M,F),d^F\right)}\right)\over b^{\cal
M}_{\left(M,F,b^F,-X\right)}}\prod_{i=0}^n\left(\det\left(\left.
D^2_{b_T}\right|_{\Omega^*_{[0,1],T}(M,F)^\perp\cap
\Omega^{i}(M,F)}\right)\right)^{(-1)^ii}\\
={P_\infty^{\det H}\left(b^{\rm
RS}_{\left(M,F,g^{TM},b^F\right)}\right)\over  b^{\cal
M}_{\left(M,F,b^F,-X\right)}} \exp\left(-2T{\rm
rk}(F)\int_Mfe\left(TM,\nabla^{TM}\right)\right).
\end{multline}

From now on, we will write $a\simeq b$ for $a,\ b\in{\bf C}$ if
$e^a=e^b$. Thus, we can rewrite (\ref{3.20}) as
\begin{multline}\label{3.21}
 \log\left({P_\infty^{\det H}\left(b^{\rm
RS}_{\left(M,F,g^{TM},b^F\right)}\right)\over  b^{\cal
M}_{\left(M,F,b^F,-X\right)}}\right)\simeq\log\left(
{P_T^{[0,1],\det H}\left(b_{\det
H^*\left(\Omega^{\ast}_{[0,1],T}(M,F),d^F\right)}\right)\over
b^{\cal M}_{\left(M,F,b^F,-X\right)}}\right)\\
+\sum_{i=0}^n(-1)^ii \log\left(\det\left(\left.
D^2_{b_T}\right|_{\Omega^{*}_{[0,1],T}(M,F)^\perp\cap
\Omega^{i}(M,F)}\right)\right)\\ +2T{\rm rk
}(F)\int_Mfe\left(TM,\nabla^{TM}\right) .
\end{multline}

Let $T_0>0$ be as in Theorem \ref{t3.4}. For any $T\geq T_0$ and
$s\in {\bf C}$ with ${\rm Re}(s)\geq n+1$, set
\begin{align}\label{3.22}
 \theta_T(s)={1\over \Gamma(s)}\int_0^{+\infty}t^{s-1}{\rm
 Tr}_s\left[N\exp\left(-tD_{b_T}^2\right)P^{(1,+\infty)}_T\right]dt.
\end{align}
By (\ref{3.13}), $\theta_T(s)$ is well-defined and can be extended
to a meromorphic function which is holomorphic at $s=0$ (cf.
\cite{S}). Moreover,
\begin{align}\label{3.23}
\sum_{i=0}^n(-1)^ii \log\left(\det\left(\left.
D^2_{b_T}\right|_{\Omega^{*}_{[0,1],T}(M,F)^\perp\cap
\Omega^{i}(M,F)}\right)\right)\simeq -\left. {\partial
\theta_T(s)\over \partial s}\right|_{s=0}.
\end{align}

Let $d=\alpha^2$ with $\alpha$ being as in Theorem \ref{t3.9}.
From (\ref{3.22}) and Theorems \ref{t3.4}-\ref{t3.6}, one finds
\begin{multline}\label{3.24}
 \left. {\partial \theta_T(s)\over \partial s}\right|_{s=0}=\int_0^d\left(
 {\rm
 Tr}_s\left[N\exp\left(-tD_{b_T}^2\right)P^{(1,+\infty)}_T\right]-{a_{-1}\over\sqrt{t}}-{n\over 2}\chi(F)+
  {\chi}'(F)\right) {dt\over t}\\
 +\int_d^{+\infty}{\rm
 Tr}_s\left[N\exp\left(-tD_{b_T}^2\right)P^{(1,+\infty)}_T\right]{dt\over
 t}
 - {2a_{-1}\over{\sqrt{d}}}\\ -\left(\Gamma'(1)-{\log d}\right)\left({n\over 2}\chi(F)-  {\chi}'(F)\right),
\end{multline}
where we denote for simplicity that
\begin{align}\label{3.25}
a_{-1}={\rm rk}(F)\int_M\int^BL\exp\left(-{\dot{R}^{TM}\over
2}\right).
\end{align}

\begin{prop}\label{t3.10} One has
\begin{align}\label{3.26}
 \lim_{T\rightarrow +\infty} \int_d^{+\infty} {\rm
 Tr}_s\left[N\exp\left(-tD_{b_T}^2\right)P^{(1,+\infty)}_T\right]{dt\over t}=0.
\end{align}
\end{prop}

{\it Proof}. This follows from Theorem \ref{t3.4} directly. \ \
Q.E.D.

$\ $

Now we write
\begin{multline}\label{3.27}
  \int_0^d\left(
 {\rm
 Tr}_s\left[N\exp\left(-tD_{b_T}^2\right)P^{(1,+\infty)}_T\right]-{a_{-1}\over\sqrt{t}}-{n\over 2}\chi(F)+
  {\chi}'(F)\right) {dt\over t} \\
 =\int_0^d\left(
 {\rm
 Tr}_s\left[N\exp\left(-tD_{b_T}^2\right)\right]-
 {a_{-1}\over\sqrt{t}}-{n\over 2}\chi(F) \right) {dt\over t}\\
 -\int_0^d\left(
 {\rm
 Tr}_s\left[N\exp\left(-tD_{b_T}^2\right)P^{[0,1]}_T\right]-
  {\chi}'(F)\right) {dt\over t}.
\end{multline}

From Theorem \ref{t3.5}, one deduces that
\begin{align}\label{3.28}
 \lim_{T\rightarrow +\infty}\int_0^d\left(
 {\rm
 Tr}_s\left[N\exp\left(-tD_{b_T}^2\right)P^{[0,1]}_T\right]-
  {\chi}'(F)\right) {dt\over t}  =0.
\end{align}

To study the first term in the right hand side  of (\ref{3.27}),
we observe first that for any $T\geq 0$, one has
\begin{align}\label{3.29}
  e^{-Tf}D_{b_T}^2e^{Tf}=\left(D_b+T\widehat{c}(\nabla f)\right)^2.
\end{align}
Thus, one has
\begin{align}\label{3.30}
{\rm
 Tr}_s\left[N\exp\left(-tD_{b_T}^2\right) \right] ={\rm
 Tr}_s\left[N\exp\left(-t\left(D_{b }+T\widehat{c}(\nabla f)\right)^2\right) \right]  .
\end{align}

By (\ref{3.30}), one writes
\begin{multline}\label{3.31}
\int_0^d\left(
 {\rm
 Tr}_s\left[N\exp\left(-tD_{b_T}^2\right)\right]-
 {a_{-1}\over\sqrt{t}}-{n\over 2}\chi(F) \right) {dt\over t}\\
 =2\int_0^{\sqrt{d}}\left(
 {\rm
 Tr}_s\left[N\exp\left(-\left(tD_{b }+tT\widehat{c}(\nabla f)\right)^2\right) \right] -
 {a_{-1}\over {t}}-{n\over 2}\chi(F)\right) {dt\over t}\\
 =2\int_{1\over \sqrt{T}}^{\sqrt{d}}\left(
 {\rm
 Tr}_s\left[N\exp\left(-\left(tD_{b }+tT\widehat{c}(\nabla f)\right)^2\right) \right] -
 {a_{-1}\over {t}}-{n\over 2}\chi(F)\right) {dt\over t}\\
 +2\int_0^{1\over \sqrt{T}}\left(
 {\rm
 Tr}_s\left[N\exp\left(-\left(tD_{b }+tT\widehat{c}(\nabla f)\right)^2\right) \right] -
 {a_{-1}\over {t}}-{n\over 2}\chi(F)\right) {dt\over t}\\
=2\int_{1}^{ \sqrt{dT}}\left(
 {\rm
 Tr}_s\left[N\exp\left(-\left({t\over \sqrt{T}}D_{b }+t\sqrt{T}\widehat{c}(\nabla f)\right)^2\right) \right] -
 {\sqrt{T}\over {t}}a_{-1}-{n\over 2}\chi(F)\right) {dt\over t}\\
 +2\int_0^{1\over \sqrt{T}}\left(
 {\rm
 Tr}_s\left[N\exp\left(-\left(tD_{b }+tT\widehat{c}(\nabla f)\right)^2\right) \right] -
 {a_{-1}\over {t}}-{n\over 2}\chi(F)\right) {dt\over t} .
\end{multline}

In view of Theorem \ref{t3.7}, we write
\begin{multline}\label{3.32}
\int_0^{1\over \sqrt{T}}\left(
 {\rm
 Tr}_s\left[N\exp\left(-\left(tD_{b }+tT\widehat{c}(\nabla f)\right)^2\right) \right] -
 {a_{-1}\over {t}}-{n\over 2}\chi(F)\right) {dt\over t}\\
 = \int_0^{1\over \sqrt{T}}\left(
 {\rm
 Tr}_s\left[N\exp\left(-\left(tD_{b }+tT\widehat{c}(\nabla f)\right)^2\right) \right] -
 {{1}\over {t}}\int_M\int^BL\exp\left(-B_{{(tT)}^2}\right){\rm rk}(F)\right.\\
 \left.-{tT\over 2}\int_M\theta\left(F,b^F\right)\int^B\widehat{df}
 \exp\left(-B_{(tT)^2}\right)-{n\over 2}\chi(F)\right) {dt\over t}\\
 + \int_0^{1\over \sqrt{T}}\left({{1}\over {t}}\int_M\int^BL\exp\left(-B_{{(tT)}^2}\right){\rm rk}(F)
 -{a_{-1}\over t}\right){dt\over t} \\
 +\int_0^{1\over \sqrt{T}}{tT\over 2}\int_M\theta\left(F,b^F\right)\int^B\widehat{df}
 \exp\left(-B_{(tT)^2}\right){dt\over t}.
\end{multline}

By \cite[Definitions 3.6, 3.12 and Theorem 3.18]{BZ1}, one has, as
$T\rightarrow +\infty$,
\begin{multline}\label{3.33}
\int_0^{1\over \sqrt{T}}{tT\over
2}\int_M\theta\left(F,b^F\right)\int^B\widehat{df}
 \exp\left(-B_{(tT)^2}\right){dt\over t}\\
 ={1\over 2}
 \int_0^{  \sqrt{T}} \int_M\theta\left(F,b^F\right)\int^B\widehat{df}
 \exp\left(-B_{t^2}\right){dt }  \\
 \rightarrow{1\over 2}
 \int_0^{  +\infty} \int_M\theta\left(F,b^F\right)\int^B\widehat{df}
 \exp\left(-B_{t^2}\right){dt }\\
 ={1\over 2}\int_M\theta\left(F,b^F\right)
 (\nabla f)^*\psi\left(TM,\nabla^{TM}\right) .
\end{multline}

By \cite[(3.58)]{BZ1}, we have for any $T\geq 0$ that
\begin{multline}\label{3.34}
 \int_M\int^B\left(L\exp\left(-B_{ T }\right)
 -L\exp\left(-B_{ 0 }\right)\right)  \\
 =-\sqrt{T}f \int_M\int^B\left( \exp\left(-B_{ T }\right)-  \exp\left(-B_{{0}}\right)
 \right)\\
 + \int_M{f\over 2}\int_0^T\left(\int^B\left(\exp\left(-B_{ t }\right)
 -\exp\left(-B_{{0}}\right)\right)\right){dt\over \sqrt{t}}.
\end{multline}

From (\ref{3.34}), one deduces easily that
\begin{align}\label{3.35}
\lim_{T\rightarrow0^+}{1\over\sqrt{T}}\int_M\int^B\left(L\exp\left(-B_{{
T } }\right)- L\exp\left(-B_{{0}}\right)
 \right)=0.
\end{align}

From \cite[(3.54)]{BZ1}, (\ref{3.35}) and the integration by
parts, we have
\begin{multline}\label{3.36}
\int_0^{1\over \sqrt{T}}\left({{1}\over
{t}}\int_M\int^BL\exp\left(-B_{{(tT)}^2}\right){\rm rk}(F)
 -{a_{-1}\over t}\right){dt\over t} \\
=-T{\rm rk}(F) \int_0^{ {T}} \int_M\int^B\left( L\exp\left(-B_{ t
}\right)
 - L\exp\left(-B_{ 0 }\right)\right)  {d{1\over \sqrt{t}}}\\
 =-\sqrt{T}{\rm rk}(F)\int_M \int^B\left(L\exp\left(-B_{ T
}\right)
 -L\exp\left(-B_{ 0 }\right)\right)\\ -{T}{\rm rk}(F)\int_M\int_0^Tf{\partial \over
 \partial t}\int^B\exp\left(-B_t\right)dt\\
 = -\sqrt{T}{\rm rk}(F)\int_M \int^B L\exp\left(-B_{ T
}\right)+\sqrt{T}a_{-1}-T{\rm
rk}(F)\int_Mf\int^B\exp\left(-B_T\right)
\\ +T{\rm rk}(F)\int_Mf\int^B\exp\left(-B_0\right).
\end{multline}

From Theorems \ref{t3.7}, \ref{t3.8}, (\ref{3.35}), \cite[Theorem
3.20]{BZ1}, \cite[(7.72) and (7.73)]{BZ1} and the dominate
convergence, one finds that as $T\rightarrow +\infty$,
\begin{multline}\label{3.37}
\int_0^{1\over \sqrt{T}}\left(
 {\rm
 Tr}_s\left[N\exp\left(-\left(tD_{b }+tT\widehat{c}(\nabla f)\right)^2\right) \right] -
 {{1}\over {t}}\int_M\int^BL\exp\left(-B_{{(tT)}^2}\right){\rm rk}(F)\right.\\
 \left.-{tT\over 2}\int_M\theta\left(F,b^F\right)\int^B\widehat{df}
 \exp\left(-B_{(tT)^2}\right)-{n\over 2}\chi(F)\right) {dt\over
 t}\\ =\int_0^{1 }\left(
 {\rm
 Tr}_s\left[N\exp\left(-\left({t\over\sqrt{T}}D_{b }+t\sqrt{T }
 \widehat{c}(\nabla f)\right)^2\right) \right]\right.\\ \left. -
 {\sqrt{T }\over {t}}\int_M\int^BL\exp\left(-B_{{(t\sqrt{T })}^2}\right){\rm rk}(F)\right.\\
 \left.-{t\sqrt{T }\over 2}\int_M\theta\left(F,b^F\right)\int^B\widehat{df}
 \exp\left(-B_{(t\sqrt{T })^2}\right)-{n\over 2}\chi(F)\right) {dt\over
 t}\\
 \rightarrow \int_0^1\left({1\over
 1-e^{-2t^2}}\left(\left(1+e^{-2t^2}\right)\chi'(F)-ne^{-2t^2}\chi(F)\right)\right.  \\
 +\left. {{\rm rk}(F)\over 2t^2}\sum_{x\in B}(-1)^{{\rm ind}(x)}\left(n-2\,{\rm
 ind}(x)\right)-{n\over 2}\chi(F)\right) {dt\over t}\\
 ={1\over 2}\left(\chi'(F)-{n\over
 2}\chi(F)\right)\int_0^1\left({1+e^{-2t }\over
 1-e^{-2t }}-{1\over t}\right){dt\over t}.
\end{multline}

On the other hand, by Theorems \ref{t3.8}, \ref{t3.9} and the
dominate convergence, we have that as $T\rightarrow +\infty$,
\begin{multline}\label{3.38}
 \int_1^{\sqrt{Td}}\left(
 {\rm
 Tr}_s\left[N\exp\left(-\left({t\over \sqrt{T}}D_{b }+t\sqrt{T}\widehat{c}(\nabla f)\right)^2\right) \right] -
 {\sqrt{T} \over {t}}a_{-1}-{n\over 2}\chi(F)\right) {dt\over t}\\
 = \int_1^{\sqrt{Td}}\left(
 {\rm
 Tr}_s\left[N\exp\left(-\left({t\over\sqrt{T}}D_{b }+t\sqrt{T}\widehat{c}(\nabla f)\right)^2\right) \right]
 - {\chi}'(F)
 \right) {dt\over t}\\
 +{1\over 2}  {\chi}'(F) \log \left(Td\right)+a_{-1}\sqrt{T}\left({1\over\sqrt{Td}}-1\right)
  -{n\over 4}\chi(F)\log \left(Td\right)\\
  =\int_1^{+\infty} \left({1\over
 1-e^{-2t^2}}\left(\left(1+e^{-2t^2}\right)\chi'(F)-ne^{-2t^2}\chi(F)\right)-\chi'(F)\right){dt\over t}\\
 +{1\over 2} \left(\chi'(F)-{n\over
 2}\chi(F)\right)\log \left(Td\right) +{a_{-1}\over{\sqrt{d}}}-\sqrt{T}a_{-1}+o(1) \\
 =\left(\chi'(F)-{n\over
 2}\chi(F)\right)\int_1^{+\infty}  {e^{-2t }\over
 1-e^{-2t }} {dt\over t}
 +{1\over 2} \left(\chi'(F)-{n\over
 2}\chi(F)\right)\log \left(Td\right)\\ +{a_{-1}\over\sqrt{d}}-\sqrt{T}a_{-1}+o(1).
\end{multline}

Combining (\ref{3.11}), (\ref{3.21}), (\ref{3.23})-(\ref{3.28}),
(\ref{3.31})-(\ref{3.33}) and (\ref{3.36})-(\ref{3.38}), one
deduces, by setting $T\rightarrow +\infty$, that
\begin{multline}\label{3.39}
\log\left({P_\infty^{\det H}\left(b^{\rm
RS}_{\left(M,F,g^{TM},b^F\right)}\right)\over  b^{\cal
M}_{\left(M,F,b^F,-X\right)}}\right)\simeq -2\,{\rm rk}(F){\rm
Tr}_s^B[f] T+\left(\chi'(F)-{n\over
 2}\chi(F)\right)\log T\\
 -\left(\chi'(F)-{n\over
 2}\chi(F)\right)\log \pi -\int_M\theta\left(F,b^F\right)(\nabla
 f)^*\psi\left(TM,\nabla^{TM}\right)\\
 +2\sqrt{T}{\rm
 rk}(F)\int_M\int^BL\exp\left(-B_T\right)-2\sqrt{T}a_{-1}+2\,T{\rm
 rk}(F)\int_Mf\int^B\exp\left(-B_T\right)\\ -2\,T{\rm
 rk}(F)\int_Mf\int^B\exp\left(-B_0\right)-\left(\chi'(F)-{n\over
 2}\chi(F)\right)\int_0^1\left({1+e^{-2t}\over 1-e^{-2t}}-{1\over
 t}\right){dt\over t}\\
 -\left(\chi'(F)-{n\over
 2}\chi(F)\right)\int_1^{+\infty} {2\, e^{-2t}\over 1-e^{-2t}} {dt\over
 t}-\left(\chi'(F)-{n\over
 2}\chi(F)\right)\log \left(Td\right)
 - {2a_{-1}\over{\sqrt{d}}}\\ +2\sqrt{T}a_{-1}+2\,T{\rm
 rk}(F)\int_Mf\, e\left(TM,\nabla^{TM}\right)+ {2a_{-1}\over{\sqrt{d}}}\\
 -\left(\Gamma'(1)- \log d\right)\left(\chi'(F)-{n\over
 2}\chi(F)\right)+o(1)\\
 = 2\,T{\rm rk}(F)\int_Mf\left(\int^B\exp\left(-B_T\right)-\sum_{x\in
 B}(-1)^{{\rm ind}(x)}\delta_x\right)\\
  -\left(\chi'(F)-{n\over
 2}\chi(F)\right)\left(\int_0^1\left({1+e^{-2t}\over 1-e^{-2t}}-{1\over
 t}\right){dt\over t}+\int_1^{+\infty} {2\, e^{-2t}\over 1-e^{-2t}} {dt\over
 t}\right)\\
 -\left(\chi'(F)
 -{n\over
 2}\chi(F)\right)\left(\log\pi+\Gamma'(1)\right)+2\sqrt{T}{\rm rk}\left(F\right)\int_M\int^BL\exp\left(-B_T\right)
 \\ -\int_M\theta\left(F,b^F\right)(\nabla
 f)^*\psi\left(TM,\nabla^{TM}\right)+o(1).
\end{multline}

By \cite[Theorem 3.20]{BZ1} and \cite[(7.72)]{BZ1}, one has
\begin{multline}\label{3.40}
\lim_{T\rightarrow +\infty} 2\,T{\rm
rk}(F)\int_Mf\left(\int^B\exp\left(-B_T\right)-\sum_{x\in
 B}(-1)^{{\rm ind}(x)}\delta_x\right)\\ =-\left(\chi'(F)
 -{n\over
 2}\chi(F)\right),
\end{multline}
\begin{align}\label{3.41}
\lim_{T\rightarrow +\infty}2\sqrt{T}{\rm
rk}(F)\int_M\int^BL\exp\left(-B_T\right)=2\left(\chi'(F)
 -{n\over
 2}\chi(F)\right).
\end{align}

On the other hand, by \cite[(7.93)]{BZ1}, one has
\begin{align}\label{3.42}
 \int_0^1\left({1+e^{-2t}\over 1-e^{-2t}}-{1\over
 t}\right){dt\over t}+\int_1^{+\infty} {2\, e^{-2t}\over 1-e^{-2t}} {dt\over
 t}  =1-\log\pi-\Gamma'(1).
\end{align}

From (\ref{3.39})-(\ref{3.42}), we get (\ref{3.7}), which
completes the proof of Theorem \ref{t3.1}.\ \ Q.E.D.

$\ $

\begin{Rem}\label{t3.11} One finds that we have used the strategy
outlined in \cite[Appendix]{BZ2} to prove Theorem \ref{t3.1},
instead of using that in \cite[Section 7]{BZ1}.  In particular, we
avoid the explicit use of \cite[Theorem 3.9]{BZ1} which is crucial
in \cite[Section 7]{BZ1}, though we still make use of the
variation formulas \cite[(3.54) and (3.58)]{BZ1}.
\end{Rem}

\begin{Rem}\label{t3.12} By Theorem \ref{t3.7}, one deduces that
\begin{multline}\label{3.43}\lim_{T\rightarrow +\infty}\int_0^{1 }\left(
 {\rm
 Tr}_s\left[N\exp\left(-\left({t\over\sqrt{T}}D_{b }+t\sqrt{T }
 \widehat{c}(\nabla f)\right)^2\right) \right]\right.\\ \left. -
 {\sqrt{T }\over {t}}\int_M\int^BL\exp\left(-B_{{(t\sqrt{T })}^2}\right){\rm rk}(F)\right.\\
 \left.-{t\sqrt{T }\over 2}\int_M\theta\left(F,b^F\right)\int^B\widehat{df}
 \exp\left(-B_{(t\sqrt{T })^2}\right)-{n\over 2}\chi(F)\right) {dt\over
 t}=0. \end{multline}
 Combining with (\ref{3.37}), one gets
\begin{align}\label{3.44}
  \int_0^1\left({1+e^{-2t}\over 1-e^{-2t}}-{1\over
 t}\right){dt\over t}=0.
\end{align}
\end{Rem}

\section{Asymptotics of the symmetric bilinear torsion of the Witten complex} \label{s4}
\setcounter{equation}{0}

 In this section, we prove Theorems 3.3 and 3.5.

 We make the same
assumptions and use the same notations as in Section \ref{s3}.

\subsection{\normalsize   Some formulas related to $D_b$}
\label{s4.1}

Recall that $b^F$ is a nondegenerate symmetric bilinear form on a
complex flat vector bundle $F$ over an oriented closed Riemannian
manifold $M$. Then it determines a nondegenerate symmetric
bilinear form $\langle\ ,\ \rangle_b$ on $\Omega^*(M,F)$ (cf.
(\ref{2.17})).

Recall that the formal adjoint $d_b^{F*}$ of $d^F$ with respect to
the symmetric bilinear form $\langle\ ,\ \rangle_b$   has been
defined in (\ref{2.19}), and $D_b$ is the operator defined by
\begin{align}\label{4.1}
  D_b=d^F+d_b^{F*}.
\end{align}

Let
\begin{align}\label{4.2}
  \omega_b^F=\omega_b\left(F,\nabla^F\right)=\left(b^F\right)^{-1}\nabla^Fb^F
\end{align}
be defined as in \cite{BH2}.

Let $\nabla=\nabla^{\Lambda^*(T^*M\otimes F)}$ be the tensor
product connection on $\Lambda^*(T^*M)\otimes F$ obtained from the
Levi-Civita connection $\nabla^{TM}$ associated to $g^{TM}$ and
the flat connection $\nabla^F$ on $F$.

For any $X\in TM$, let $X^*\in T^*M$ corresponds to $X$ via
$g^{TM}$. Recall that
\begin{align}\label{4.3}
 c(X)=X^*-i_X,\ \ \ \ \ \widehat{c}(X)=X^*+i_X
\end{align}
denote the Clifford actions on $\Lambda^*(T^*M )$, where $X^*$ and
$i_X$ are the exterior and interior multiplications respectively
(cf. \cite[Section 4]{BZ1}).

For any oriented orthonormal basis $e_1,\, \dots,\, e_n$ of $TM$,
set
\begin{align}\label{4.4}
  c\left(\omega_b^F\right)=\sum_{i=1}^n
  c\left(e_i\right)\omega_b^F\left(e_i\right),\ \ \ \ \ \widehat{c}\left(\omega_b^F\right)=\sum_{i=1}^n
  \widehat{c}\left(e_i\right)\omega_b^F\left(e_i\right).
\end{align}

With these definitions and notations one verifies easily that (cf.
\cite[(92)]{BH2})
\begin{align}\label{4.5}
d^F+d_b^{F*}=\sum_{i=1}^n c\left(e_i\right)\nabla_{e_i}+{1\over
2}c\left(\omega_b^F\right)-{1\over
2}\widehat{c}\left(\omega_b^F\right).
\end{align}

Recall that $g^F$ is a Hermitian metric on $F$. Together with
$g^{TM}$ it determines an inner product $\langle\ ,\ \rangle_g$ on
$\Omega^*(M,F)$ (cf. \cite[(2.2)]{BZ1} and \cite[(2.3)]{BZ2}).

Let $d_g^{F*}$ be the formal adjoint of $d^F$ with respect to
$\langle\ ,\ \rangle_g$.

Set as in \cite{BZ1} and \cite{BZ2} that
\begin{align}\label{4.6}
\omega_g^F=\omega_g\left(F,\nabla^F\right)=\left(g^F\right)^{-1}\nabla^Fg^F.
\end{align}
Then $\omega_g^F$ is a one form taking values in the self-adjoint
elements in ${\rm End}(F)$. Moreover,
\begin{align}\label{4.7}
 \nabla^{F,u}=\nabla^F+{1\over 2}\omega_g^F
\end{align}
is a Hermitian connection on $F$ with respect to $g^F$ (cf.
\cite[Section 4]{BZ1} and \cite[Section 2]{BZ2}). Let $\nabla^{
u}$ be the associated tensor product connection on
$\Lambda^*(T^*M)\otimes F$.

By \cite[(4.25)]{BZ1}, one has
\begin{multline}\label{4.8}
  D_g:=d^F+d_g^{F*}= \sum_{i=1}^n
  c\left(e_i\right)\nabla_{e_i}^u-{1\over 2}\widehat{c}\left(\omega_g^F\right)\\ = \sum_{i=1}^n
  c\left(e_i\right)\nabla_{e_i}+{1\over
2}c\left(\omega_g^F\right)-{1\over
2}\widehat{c}\left(\omega_g^F\right).
\end{multline}

From (\ref{4.5}) and (\ref{4.8}), one gets
\begin{align}\label{4.9}
 d^F+d^{F*}_b=d^F+d^{F*}_g+{1\over
2}c\left(\omega_b^F\right)-{1\over
2}\widehat{c}\left(\omega_b^F\right)-{1\over
2}c\left(\omega_g^F\right)+{1\over
2}\widehat{c}\left(\omega_g^F\right) .
\end{align}

Write $\omega^F_b$ as
\begin{align}\label{4.10}
 \omega^F_b=\omega^F_{b,1}+\omega^F_{b,2}  ,
\end{align}
where $\omega^F_{b,1}$ (resp.  $\omega^F_{b,2}$) takes values in
self-adjoint (resp. skew-adjoint) elements (with respect to $g^F$)
in ${\rm End}(F)$.

From (\ref{4.9}), one gets the decomposition of $D_b$ into
self-adjoint and skew-adjoint parts (with respect to $\langle\ ,\
\rangle_g$) as follows,
\begin{multline}\label{4.11}
d^F+d_b^{F*}=  \left( d^F+d_g^{F*}+{1\over
2}\widehat{c}\left(\omega_g^F\right)-{1\over
2}\widehat{c}\left(\omega_{b,1}^F\right)
   +{1\over 2}c\left(\omega_{b,2}^F\right)\right)\\ +
  \left( -{1\over
2}c\left(\omega_g^F\right) +{1\over
2}c\left(\omega_{b,1}^F\right)-{1\over
2}\widehat{c}\left(\omega_{b,2}^F\right)\right).
\end{multline}

\subsection{\normalsize   Witten deformation and some basic estimates}
\label{s4.2}

Let $f:M\rightarrow {\bf R}$ be a Morse function on $M$. We make
the assumption that the Riemannian metric $g^{TM}$ and $f$ verify
the condition (\ref{3.1}). We also assume that $g^F$, like $b^F$,
is flat near the set of critical points of $f$.

Following Witten \cite{W}, for any $T\in{\bf R}$, set
\begin{align}\label{4.12}
 d_T^F=e^{-Tf}d^Fe^{Tf},\ \ \ \delta_{b,T}^F=e^{Tf}d^{F*}_be^{-Tf},\ \ \
 \delta_{g,T}^F=e^{Tf}d^{F*}_ge^{-Tf}.
\end{align}
Set
\begin{align}\label{4.13}
 \widetilde{D}_{b,T}=d_T^F+\delta_{b,T}^F= D_b+T\widehat{c}(df), \ \
 \  \widetilde{D}_{g,T}=d_T^F+\delta_{g,T}^F= D_g+T\widehat{c}(df)
  .
\end{align}
Observe that the skew-adjoint part of $\widetilde{D}_{b,T}$ is the
same as that of $\widetilde{D}_{b }$.

Let $\|\ \|_0$ be the $L^2$ norm on $\Omega^*(M,F)$ associated to
$\langle\ ,\ \rangle_g$. For any $q>0$, let $\|\ \|_q$ be a fixed
$q$-Sobolev norm on $\Omega^*(M,F)$.

\begin{prop} \label{t4.1} For any open neighborhood $U$ of $B$,
there exist $T_0>0$, $C>0$, $c>0$ such that for any $s\in
\Omega^*(M,F)$ with ${\rm supp}(s)\subset M\setminus U$  and
$T\geq T_0$, one has
\begin{align}\label{4.14}
 \left\| \widetilde{D}_{b,T}s\right\|_0^2\geq
 C\left(\|s\|_1^2+(T-c)\|s\|_0^2\right)
  .
\end{align}
\end{prop}

{\it Proof.} From (\ref{4.11}) and (\ref{4.13}), one sees that the
formal adjoint $\widetilde{D}_{b,T}^*$ of $\widetilde{D}_{b,T}$ is
given by
\begin{multline}\label{4.15}
\widetilde{D}_{b,T}^* =  \left( D_{g }+T\widehat{c}(df)+{1\over
2}\widehat{c}\left(\omega_g^F\right)-{1\over
2}\widehat{c}\left(\omega_{b,1}^F\right)
   +{1\over 2}c\left(\omega_{b,2}^F\right)\right)\\ -
  \left( -{1\over
2}c\left(\omega_g^F\right) +{1\over
2}c\left(\omega_{b,1}^F\right)-{1\over
2}\widehat{c}\left(\omega_{b,2}^F\right)\right). \end{multline}

For simplicity, we denote by
\begin{align}\label{4.16}
  A^F={1\over
2}\widehat{c}\left(\omega_g^F\right)-{1\over
2}\widehat{c}\left(\omega_{b,1}^F\right)
   +{1\over 2}c\left(\omega_{b,2}^F\right),
\end{align}
$$B^F=-{1\over
2}c\left(\omega_g^F\right) +{1\over
2}c\left(\omega_{b,1}^F\right)-{1\over
2}\widehat{c}\left(\omega_{b,2}^F\right).$$ Then one computes
\begin{multline}\label{4.17}
\widetilde{D}_{b,T}^*\widetilde{D}_{b,T}= \left(D_{g
}+A^F\right)^2+\left(D_{g }+A^F\right)B^F-B^F\left(D_{g
}+A^F\right)-\left(B^F\right)^2 \\ +T\left(\left[ D_{g }+A^F
,\widehat{c}(df)\right]+\widehat{c}(df)B^F-B^F\widehat{c}(df)\right)+T^2|df|^2,
 \end{multline}
 where by $[\ , \ ]$ we denote the super bracket in the sense of
 Quillen \cite{Q1}.

 Since it is easy to check (cf. \cite[(5.17)]{BZ1}) that
\begin{align}\label{4.18}
\left[ D_{g }
,\widehat{c}(df)\right]=\sum_{i=1}^nc\left(e_i\right)\widehat{c}\left(\nabla^{TM}_{e_i}\nabla
f\right)-\omega_g^F(\nabla f),
\end{align}
where $\nabla f\in \Gamma(TM)$ is the gradient vector field of $f$
with respect to $g^{TM}$, is of order zero, the coefficient of $T$
in the right hand side of (\ref{4.17}) is of order zero.

Also, it is clear that there is $c_0>0$ such that for any $x\in
M\setminus U$,
\begin{align}\label{4.19}
 |df(x)|\geq c_0.
\end{align}

From (\ref{4.17}) and (\ref{4.19}), one gets Proposition
\ref{t4.1} easily, as
\begin{align}\label{4.20}
 \left\| \widetilde{D}_{b,T}s\right\|_0^2=\left\langle
  \widetilde{D}_{b,T}s,\widetilde{D}_{b,T} s\right\rangle
  =\left\langle
  \widetilde{D}_{b,T}^*\widetilde{D}_{b,T}s,  s\right\rangle.
\end{align}
Q.E.D.

\begin{prop} \label{t4.2} For any $c >0$, there exists $T_c>0$ such
that for any $T\geq T_c$, $z\in {\bf C}$ with $|z|=c$, $z\not\in
{\rm Spec}(\widetilde{D}_{b,T}^2)$.
\end{prop}

{\it Proof.} For any $p\in B$, let ${\bf y}=(y_1,\dots,y_n)$ be
the coordinate system  of $p$ as in (\ref{3.1}), in an open ball
$U_p$ of radius $4a$, around $p$. We also assume that both $b^F$
and $g^F$ are flat on each  $U_p$, $p\in B$. The existence of
$a>0$ is clear.

By (\ref{4.9}), one then has
\begin{align}\label{4.21}
  D_b=D_g\ \ \ {\rm on}\ \ \ U_B=\bigcup_{p\in B}U_p\ \ .
\end{align}

Let $\gamma:{\bf R}\rightarrow [0,1]$ be a smooth function such
that $\gamma(x)=1$ if $|x|\leq a$, while $\gamma(x)=0$ if $|x|\geq
2a$.

For any $T>0$ and $p\in B$, set
\begin{align}\label{4.22}
 \alpha_{p,T}=\int_{U_p} \gamma(|{\bf y}|)^2\exp\left(-T|{\bf
 y}|^2\right)dy^1\wedge\cdots \wedge dy^n,
\end{align}
$$\rho_{p,T}={\gamma(|{\bf y}|)\over
\sqrt{\alpha_{p,T}}}\exp\left(-{T|{\bf y}|^2\over 2}\right)dy^1
\wedge\cdots \wedge dy^{n_f(p)} ,$$ where $n_f(p)={\rm ind}(p)$ is
the Morse index of $f$ at $p$. Then $\rho_{p,T}\in
\Omega^{n_f(p)}(M )$ is of unit length with compact support
contained in $U_p$.

Set
\begin{align}\label{4.23}
 E_T=\bigoplus_{p\in B}\left\{\rho_{p,T}\otimes h_p: p\in B, \ \ h_p\in F_p\right\}.
 \end{align}
 Let $E_T^\perp$ be the orthogonal complement to $E_T$ in
 $L^2(\Omega^*(M,F))$ with respect to $\langle\ ,\ \rangle_g$.
 Then one has the orthogonal decomposition
\begin{align}\label{4.24}
 L^2\left(\Omega^*(M,F)\right)=E_T\oplus E_T^\perp.
 \end{align}
 Let $p_T$, $p_T^\perp$ be the orthogonal projections from
 $L^2(\Omega^*(M,F))$ onto $E_T$, $E_T^\perp$ respectively.

 Following \cite[Section 9b)]{BL} (cf. \cite[(5.19)]{Z}), set
\begin{align}\label{4.25}
 \widetilde{D}_{b,T,1}=p_T\widetilde{D}_{b,T}p_T,\ \ \
 \widetilde{D}_{b,T,2}=p_T\widetilde{D}_{b,T}p_T^\perp,
 \end{align}
$$ \widetilde{D}_{b,T,3}=p_T^\perp\widetilde{D}_{b,T}p_T,\ \ \
 \widetilde{D}_{b,T,2}=p_T^\perp\widetilde{D}_{b,T}p_T^\perp.$$

 From (\ref{4.17}), (\ref{4.20}), (\ref{4.21}), (\ref{4.25}) and
 proceed as in \cite[Section 9]{BL} and \cite[Proof of Proposition
 5.6]{Z}, one can prove in the same way  that there exist $T_0>0$, $C>0$ such that for
 any $T\geq T_0$, one has
\begin{align}\label{4.26}
\widetilde{D}_{b,T,1}=0,
 \end{align}
\begin{align}\label{4.27}\left\|\widetilde{D}_{b,T,2}s\right\|_0\leq {\|s\|_0\over T},\ \
\
 \  \left\|\widetilde{D}_{b,T,3}s'\right\|_0\leq {\|s'\|_0\over T}\end{align}
for  any $s\in E_T^\perp\cap{\bf H}^1(M,F)$,  $s'\in E_T $, where
${\bf H}^1(M,F)$ is the  Sobolev space with respect to
 the Sobolev norm $\|\ \|_1$ on $\Omega^*(M,F)$, and
\begin{align}\label{4.28}\left\|\widetilde{D}_{b,T,4}s\right\|_0\geq C\sqrt{T}{\|s\|_0 } \end{align}
for  any $ s\in E_T^\perp\cap{\bf H}^1(M,F) $.

Now for any $\lambda\in{\bf C}$, $T\geq T_0$ and
$s\in\Omega^*(M,F)$, by (\ref{4.26})-(\ref{4.28}), we have (cf.
\cite[(5.26)]{Z})
\begin{multline}\label{4.29}
 \left\|\left(\lambda-\widetilde{D}_{b,T}\right)s\right\|_0\geq {1\over
 2}\left\|\lambda p_Ts-\widetilde{D}_{b,T,2}p_T^\perp s\right\|_0+{1\over
 2}\left\|\lambda p_T^\perp s
 -\widetilde{D}_{b,T,3}s-\widetilde{D}_{b,T,4}p_T^\perp
 s\right\|_0\\
 \geq{1\over 2}\left(\left(|\lambda|-{1\over
 T}\right)\left\|p_Ts\right\|_0+\left(C\sqrt{T}-|\lambda|-{1\over
 T}\right) \left\|p_T^\perp s\right\|_0 \right).
 \end{multline}

 From (\ref{4.29}), one sees easily that there exist $C_0>0$, $T_0'\geq T_0$ such that for
 any $T\geq T_0'$ and    $\lambda\in{\bf
 C}$ with $|\lambda|^2=c$, one has
\begin{align}\label{4.30}
\left\|\left(\lambda^2-\widetilde{D}_{b,T}^2\right)s\right\|_0=
\left\|\left(\lambda+\widetilde{D}_{b,T}\right)\left(\lambda-\widetilde{D}_{b,T}\right)s\right\|_0
\geq C_0\|s\|_0,
 \end{align}
 from which Proposition \ref{t4.2} follows. \ \ Q.E.D.

$\ $

 From now on, we take $c=1$, $T_{c=1}$ as in Proposition
 \ref{t4.2} and assume $T\geq T_1$.

Let $\widetilde{\Omega}^*_{[0,1], T}(M,F)$ be defined as in
(\ref{2.23}) with respect to $\widetilde{D}_{b,T}$. Let
$\widetilde{P}^{[0,1]}_T$ be the orthogonal projection from
$L^2(\Omega^*(M,F))$ onto  $\widetilde{\Omega}^*_{[0,1], T}(M,F)$.

For any $p\in B$, let $[W^u(p)]^*$ admit a Hermitian metric such
that $|W^u(p)^*|=1$. Let $[W^u(p)]^*\otimes F_p$ carry  the tensor
product metric from the above one with $g^{F_p}$. Let $C^*(W^u,F)$
carry  a Hermitian metric through the orthogonal direct sum of the
Hermitian metrics on $[W^u(p)]^*\otimes F_p$'s.

Let $J_T:C^*(W^u,F)\rightarrow \Omega^*(M,F)$ be the isometry
defined by that for any $p\in B$, $h\in F_p$ and ${\bf y}$ the
coordinate system as above in $U_p$,
\begin{align}\label{4.31}
J_T\left(W^u(p)^*\otimes h\right)({\bf y})= \rho_{p,T}\otimes h.
 \end{align}

 From (\ref{4.11}) and (\ref{4.21}), one can proceed in exactly
 that same way as in \cite[Theorem 8.8]{BZ1} and \cite[Theorem
 6.7]{BZ2} to get the following result.

 \begin{thm}\label{t4.3} There exists $c>0$ such that as
 $T\rightarrow +\infty$, for any $s\in C^*(W^u,F)$,
\begin{align}\label{4.32}
\left(\widetilde{P}^{[0,1]}_TJ_T-J_T\right)s=O\left(e^{-cT}\right)s\
\ \ {\rm uniformly\ on}\ M.
 \end{align}
 \end{thm}

\subsection{\normalsize  Proof of Theorems 3.5}
\label{s4.3}

From Theorem \ref{t4.3}, one gets immediately that
\begin{align}\label{4.33}
\dim \widetilde{\Omega}^*_{[0,1], T}(M,F)\geq \# B.
 \end{align}

 By (\ref{4.21}) and proceed as  in \cite[Proof of Proposition
 5.5]{Z}, one sees that indeed, (\ref{4.33}) holds in equality.

 Since $\widetilde{P}^{[0,1]}_T$ preserves the ${\bf Z}$-grading
 of $\Omega^*(M,F)$ (as $\widetilde{D}_{b,T}^2$ does), by applying
 (\ref{4.32}) in each grade and by (\ref{4.33}) with equality, one
 then gets that for any $0\leq i\leq n$,
\begin{align}\label{4.34}
\dim \widetilde{\Omega}^i_{[0,1], T}(M,F)={\rm rk}(F)M_i= {\rm
rk}(F)\cdot\,\#\left\{p\in B:{\rm ind}(p)=i\right\}.
 \end{align}

On the other hand, since the number $c$ in Proposition \ref{t4.2}
can be chosen arbitrarily small, one sees that when $T\rightarrow
+\infty$, one has
\begin{align}\label{4.35}
{\rm
Tr}\left[\widetilde{D}_{b,T}^2\widetilde{P}_T^{[0,1]}\right]\rightarrow
0.
\end{align}

Now consider the isomorphism $r_T:\Omega^*(M,F)\rightarrow
\Omega^*(M,F)$ defined by $r_T(s)=  e^{Tf}s$. Then it induces a
map preserving the corresponding symmetric bilinear forms, as well
as the inner products,
\begin{align}\label{4.36}
r_T: \left(\Omega^*(M,F), \langle\ ,\
\rangle_{b}\right)\longmapsto \left(\Omega^*(M,F), \langle\ ,\
\rangle_{b_T}\right).
\end{align}
$$r_T: \left(\Omega^*(M,F), \langle\ ,\
\rangle_{g}\right)\longmapsto \left(\Omega^*(M,F), \langle\ ,\
\rangle_{g_T}\right),$$ with $\langle\ ,\ \rangle_{g_T}$ obtained
from $g^{TM}$ and $g^F_T=e^{-2Tf}g^F$ (cf. \cite[(5.1)]{BZ1}).
 Moreover, one verifies directly that
\begin{align}\label{4.37}
r_T \widetilde{D}_{b,T}= {D}_{b_T}r_T.
\end{align}

From (\ref{4.34})-(\ref{4.37}), one gets Theorem \ref{t3.5}
immediately. \ \ Q.E.D.

\subsection{\normalsize  Proof of Theorems 3.3}
\label{s4.4}

We still assume that $T\geq T_{c=1}$, where $T_{c=1}$ verifies
Proposition \ref{t4.2}.

Let $e_T:C^*(W^u,F)\rightarrow \Omega^*_{[0,1],T}(M,F)$ be defined
by
\begin{align}\label{4.38}
e_T=r_T \widetilde{P}^{[0,1]}_TJ_T.
\end{align}

Recall that $C^*(W^u,F)$ carries a symmetric bilinear form
determined in (\ref{2.13}) and (\ref{2.14}), while
$\Omega^*_{[0,1],T}(M,F)$    carries the induced symmetric
bilinear form $\langle\ ,\ \rangle_{b_T}$. Let $e_T^\#$ be the
adjoint of $e_T$ with respect to these two symmetric bilinear
forms.

\begin{prop}\label{t4.4} There exists $c>0$ such that as
$T\rightarrow +\infty$,
\begin{align}\label{4.39}
e_T^\#e_T= 1 +O\left(e^{-cT}\right).
\end{align} In particular, when $T>0$ is large enough,
$e_T:C^*(W^u,F)\rightarrow \Omega^*_{[0,1],T}(M,F)$ is a ${\bf
Z}$-graded isomorphism.
\end{prop}

{\it Proof}. By the definition of $e_T$ and $e_T^\#$, one has that
for any $s,\ s'\in C^*(W^u,F)$,
\begin{align}\label{4.40}
 \left\langle e_T^\#e_Ts, s'\right\rangle_{b}=
 \left\langle e_Ts,e_Ts'\right\rangle_{b_T}
 =\left\langle \widetilde{P}^{[0,1]}_TJ_T
 s,\widetilde{P}^{[0,1]}_TJ_Ts'\right\rangle_{b}.
\end{align}

On the other hand,  from (\ref{4.22}) and (\ref{4.31}), one sees
directly that
\begin{align}\label{4.41}
    \left\langle  J_T
 s, J_Ts'\right\rangle_{b}
 =  \left\langle s ,s' \right\rangle_{b }.
\end{align}

From Theorem \ref{t4.3}, (\ref{4.40}), and (\ref{4.41}), one gets
(\ref{4.39}).

From Theorem \ref{t3.5}   and (4.39), one sees that when $T>0$ is
large enough, $e_T$ is an isomorphism.\ \
  Q.E.D.

$\ $

Recall that the quasi-isomorphism
$P_\infty:(\Omega^*(M,F),d^F)\rightarrow (C^*(W^u,F),\partial)$
has been defined in (\ref{3.2}). Let
$P_{\infty,T}:\Omega^*_{[0,1],T}(M,F)\rightarrow C^*(W^u,F)$ be
the restriction of $P_\infty$ on $\Omega^*_{[0,1],T}(M,F)$.

By (\ref{3.3}), one has
\begin{align}\label{4.42}
   \partial P_{\infty,T}=P_{\infty,T}d^F.
\end{align}

By Theorem \ref{t4.3} and (\ref{4.42}), one can proceed in exactly
the same way as in \cite[Proof of Theorem 6.11]{BZ2} (cf.
\cite[Section 6.4]{Z}), to get the following analogue of
\cite[Theorem 6.11]{BZ2}.

\begin{prop}\label{t4.5} There exists $c>0$ such that as
$T\rightarrow +\infty$, one has
\begin{align}\label{4.43}
 P_{\infty,T}e_T=e^{T{\cal F}}
 \left({\pi\over T}\right)^{N/2-n/4}\left(1+O\left(e^{-cT}\right)\right),
\end{align}
where ${\cal F}$ acts on $W^u(p)\otimes F_p$ with $p\in B$ by
multiplication by $f(p)$, and $N$ is the number operator acting on
$W^u(p)\otimes F_p$ with $p\in B$ by multiplication by ${\rm
ind}(p)$. In particular, for $T>0$   large enough,
$P_{\infty,T}e_T\in {\rm End}(C^*(W^u,F))$ is one to one.
\end{prop}

$\ $

From (\ref{4.42}) and Propositions \ref{t4.4}, \ref{t4.5}, one
sees that when $T>0$ is large enough,
\begin{align}\label{4.44}
 P_{\infty,T} :\left(\Omega^*_{[0,1],T}(M,F),d^F\right)\rightarrow
 \left(C^*\left(W^u,F\right),\partial\right)
\end{align} is a cochain isomorphism.

From Proposition \ref{t2.5} and (\ref{4.44}), one finds
\begin{align}\label{4.45}
 {P_T^{[0,1],\det H}\left(b_{\det
H^*(\Omega^*_{[0,1],T} (M,F),d^F)}\right)\over b^{\cal M}_{(
{M},F,b^F,-X)}  }=\prod_{i=0}^n
\det\left(\left.P_{\infty,T}^\#P_{\infty,T}\right|_{\Omega^i_{[0,1],T}(M,F)}\right)^{(-1)^{i+1}},
\end{align}
where $P_{\infty,T}^\#$ is the adjoint of $P_{\infty,T}$ with
respect to the symmetric bilinear forms $\langle\ ,\ \rangle_b$.

From Propositions \ref{t4.4} and \ref{t4.5}, one deduces that as
$T\rightarrow +\infty$,
\begin{multline}\label{4.46}
\det\left(\left.P_{\infty,T}^\#P_{\infty,T}\right|_{\Omega^i_{[0,1],T}(M,F)}\right)
\\
=\det\left(\left.e_Te_T^\#P_{\infty,T}^\#P_{\infty,T}\right|_{\Omega^i_{[0,1],T}(M,F)}\right)\cdot
{\det}^{-1}\left(\left.e_Te_T^\#
\right|_{\Omega^i_{[0,1],T}(M,F)}\right)\\
=\det\left(\left.
\left(P_{\infty,T}e_T\right)^\#P_{\infty,T}e_T\right|_{C^i(W^u,F)}\right)\cdot
{\det}^{-1}\left(\left.e_T^\#e_T \right|_{C^i(W^u,F)}\right)\\
=\det\left(\left.\left(1+O\left(e^{-cT}\right)\right)^\#\left({\pi\over
T}\right)^{N-n/2}e^{2T{\cal
F}}\left(1+O\left(e^{-cT}\right)\right)
 \right|_{C^i(W^u,F)}\right)\\ \cdot
{\det}^{-1}\left(\left.\left(1+O\left(e^{-cT}\right)\right)\right|_{C^i(W^u,F)}\right).
\end{multline}

From (\ref{4.45}) and (\ref{4.46}), one gets (\ref{3.11})
immediately.

The proof of Theorem \ref{t3.3} is completed. \ \ Q.E.D.

\section{Proof of Theorems \ref{t3.4}} \label{s5}
\setcounter{equation}{0}

In this section we prove Theorem \ref{t3.4}.

In view of (\ref{4.36}),  may restate Theorem \ref{t3.4} as
follows.

\begin{thm}\label{t5.1}   For any $t>0$,
\begin{align}\label{5.1}
\lim_{T\rightarrow +\infty}  {\rm Tr}_s\left[
N\exp\left(-t\widetilde{D}_{b,T}^2\right)\widetilde{P}_T^{(1,+\infty)}\right]=0,\end{align}
where   $\widetilde{P}_T^{(1,+\infty)}={\rm
Id}-\widetilde{P}_T^{[0,1]}$. Moreover, for any $d>0$, there exist
 $c>0$, $C>0$ and $T_0\geq 1$ such that for any $t\geq d$
and $T\geq T_0$,
\begin{align}\label{5.2}
 \left|{\rm Tr}_s\left[N\exp\left(-t\widetilde{D}_{b,T}^2\right)
 \widetilde{P}_T^{(1,+\infty)}\right]\right|\leq
 c\exp(-Ct).
 \end{align}
\end{thm}

$\ $

Set
\begin{align}\label{5.3}
c_{b,g}=1+2 \,  \max_{x\in M}\left\{\left|\left( -{1\over
2}c\left(\omega_g^F\right) +{1\over
2}c\left(\omega_{b,1}^F\right)-{1\over
2}\widehat{c}\left(\omega_{b,2}^F\right)\right)(x)\right|\right\}
.
 \end{align}

By the decomposition formula (\ref{4.11}) and by (\ref{4.13}), one
sees that for any $\lambda\in{\bf C}$ with $|{\rm
Im}(\lambda)|=c_{b,g}$, $\lambda-\widetilde{D}_{b,T}$ is
invertible.

Let $\Gamma=\Gamma_1\cup\Gamma_2$ be the union of two contours
defined by
$$\Gamma_1=\left\{ x\pm\sqrt{-1}c_{b,g}:\ 2\leq x\leq
+\infty\right\}\cup \left\{  2+\sqrt{-1}y:\ -c_{b,g}\leq y\leq
c_{b,g}\right\} ,$$
$$\Gamma_2=\left\{ x\pm\sqrt{-1}c_{b,g}:\ -\infty\leq x\leq
-2\right\}\cup \left\{ - 2+\sqrt{-1}y:\ -c_{b,g}\leq y\leq
c_{b,g}\right\} .$$ We orient $\Gamma$  anti-clockwise.

By Proposition 4.2, one sees that there exists $T_0>0$ such that
for any $T\geq T_0$,
\begin{align}\label{5.4}
  {\rm Tr}_s\left[N\exp\left(-t\widetilde{D}_{b,T}^2\right)
 \widetilde{P}_T^{(1,+\infty)}\right]={1\over 2\pi\sqrt{-1}}{\rm
 Tr}_s\left[N
 \int_\Gamma {e^{-t\lambda^2}\over \lambda-\widetilde{D}_{b,T}}d\lambda\right].
 \end{align}

 Let $C>0$ be the constant verifying (\ref{4.28}).
 Following \cite[(9.113)]{BL}, set for any $T\geq 1$
 that
\begin{align}\label{5.5}
  U_T=\left\{\lambda\in{\bf C}:1\leq |\lambda|\leq {C\sqrt{T}\over
  4}\right\}.
 \end{align}

 From (\ref{4.26})-(\ref{4.28}), (\ref{5.4}) and (\ref{5.5}), one
 can proceed as in \cite[Section 9e)]{BL} to show that there
 exists $T_1\geq T_0$ such that for any $T\geq T_1$, $\lambda\in
 U_T$, $\lambda-\widetilde{D}_{b,T}$ is invertible.
Moreover, for any integer $p\geq n+2$, there exists $C'>0$ such
that if $T\geq T_1$, $\lambda\in U_T$, the following analogue of
\cite[(9.142)]{BL} holds,
\begin{align}\label{5.6}
 \left|{\rm
 Tr}_s\left[N\left(\lambda-\widetilde{D}_{b,T}\right)^{-p}\right]
 -\lambda^{-p}\chi'(F)\right|\leq {C'\over\sqrt{T}}(1+|\lambda|)^{p+1}.
 \end{align}

 From (\ref{5.6}), one can proceed as in \cite[Sections 9g), 9h)]{BL},
 with an obvious modification, to complete the proof of Theorem
 \ref{t5.1}. \ \ Q.E.D.

\section{Proof of Theorem \ref{t3.6}} \label{s6}
\setcounter{equation}{0}

In this section, we provide a proof of Theorem \ref{t3.6}, which
computes the asymptotics of ${\rm Tr}_s[N\exp(-tD_{b_T}^2)]$ for
fixed $T\geq 0$ as $t\rightarrow 0$.

Since $T\geq 0$ is fixed, we may  well assume that $T=0$.

One way to prove Theorem \ref{t3.6} is to apply the method
developed in \cite[Sections 7 and 8]{BH2}, which deals directly
with the operator $D_b^2$. Here we will prove it as an application
of the corresponding result for $D_g^2$ established in
\cite[Theorem 7.10]{BZ1}. The basic idea is very simple: we use
Duhamel principle to express the heat operator of $D_b^2$ by using
the heat operator of $D_g^2$, then one can apply the results for
$D_g^2$ to obtain the required results for $D_b^2$ (Indeed, this
idea will also be used in later sections for other local index
estimates as well).

Set
\begin{align}\label{6.1}
  \omega^F=\omega_g^F-\omega_b^F.
\end{align}

From  (\ref{6.1}), one can rewrite (\ref{4.9}) as
\begin{align}\label{6.2}
 d^F+d^{F*}_b=d^F+d^{F*}_g+{1\over
2}\widehat{c}\left(\omega^F\right)-{1\over
2}c\left(\omega^F\right).
\end{align}

From (\ref{6.2}), one sees that
\begin{align}\label{6.3}
B_{b,g}:=D_b^2-D_g^2=\left(
d^F+d^{F*}_b\right)^2-\left(d^F+d^{F*}_g\right)^2
\end{align}
is a differential operator of first order.

By Duhamel principle, one deduces that for any $t>0$,
\begin{multline}\label{6.4}
 e^{-tD_b^2}=e^{-tD_g^2}+\sum_{k=1}^{n}(-1)^kt^k\int_{\Delta_k}e^{-t_1tD_g^2}B_{b,g}
e^{-t_2tD_g^2}\cdots B_{b,g}e^{-t_{k+1}tD_g^2}dt_1\cdots dt_k\\
+ (-1)^{n+1}t^{n+1}\int_{\Delta_{n+1}}e^{-t_1tD_g^2}B_{b,g}
e^{-t_2tD_g^2}\cdots B_{b,g}e^{-t_{n+2}tD_b^2}dt_1\cdots dt_{n+1},
\end{multline}
where $\Delta_k$, $1\leq k\leq n+1$, is the $k$-simplex defined by
$t_1+\cdots+t_{k+1}=1$, $t_1\geq 0$, $\cdots, $ $t_{k+1}\geq 0$.

\begin{prop} \label{t6.1}
As $t\rightarrow 0^+$, one has
\begin{align}\label{6.5}
t^{n+1}\int_{\Delta_{n+1}}{\rm Tr}_s\left[Ne^{-t_1tD_g^2}B_{b,g}
e^{-t_2tD_g^2}\cdots B_{b,g}e^{-t_{n+2}tD_b^2}\right]dt_1\cdots
dt_{n+1}\rightarrow 0.
\end{align}
\end{prop}

{\it Proof}. For any $r>0$, let $\|\ \|_r$ denote the Schatten
norm defined for any linear operator $A$ by
\begin{align}\label{6.6}
\|A\|_r=\left({\rm Tr}\left[\left(A^*A\right)^{r\over
2}\right]\right)^{1\over r}.
\end{align}

Recall the basic properties of $\|\ \|_r$ (cf. \cite{Si}) that

(i) If $A$ is of trace class, then
\begin{align}\label{6.7}
|{\rm Tr}[A]|\leq \|A\|_1,\ \ \ \ \|A\|\leq \|A\|_1.
\end{align}

(ii) For any $r>0$ and compact operator $A$ and any bounded
operator $B$,
\begin{align}\label{6.8}
\|AB\|_r\leq \|B\|\, \|A\|_r,\ \ \ \ \|BA\|_r\leq\|B\|\,\| A\|_r.
\end{align}

(iii) (H\"older inequality) For any $p,\ q,\ r>0$ with ${1\over
r}={1\over p}+{1\over q}$,
\begin{align}\label{6.9}
\|AB\|_r\leq   \|A\|_p \|B\|_q .
\end{align}

\begin{lemma}\label{t6.2} For any $r>0$,  one has as $t\rightarrow
0^+$ that
\begin{align}\label{6.10}
\left\| \exp\left(-tD_b^2\right)\right\|_r=O\left({1\over
t^{n\over 2}}\right) .
\end{align}
\end{lemma}

{\it Proof}. Since $B_{b,g}$ is of order one, by \cite[Lemma
2.8]{CH} and \cite[Lemma 1]{F}, there exists a (fixed) constant
$C>0$ such that for any $u>0$, $t>0$ with $ut\leq 1$,
\begin{align}\label{6.11}
\left\| \exp\left(-utD_g^2\right)B_{b,g}\right\|_{u^{-1}}\leq
 C (ut)^{-{1\over 2}}\left({\rm Tr}\left[\exp\left(-{tD_g^2\over 2}\right)\right]\right)^u.
\end{align}

From (\ref{6.8}), (\ref{6.9}) and (\ref{6.11}), one sees that for
any $k\geq 1$ and $(t_1,\cdots, t_{k+1})\in\Delta_k\setminus
\{t_1\cdots t_{k+1}=0\}$,
\begin{multline}\label{6.12}
  \left\|e^{-t_1tD_g^2}B_{b,g}
e^{-t_2tD_g^2}\cdots B_{b,g}e^{-t_{k+1}tD_g^2}\right\|_1\\
\leq \left\|e^{-t_1tD_g^2}B_{b,g}\right\|_{t_1^{-1}}\cdots
\left\|e^{-t_ktD_g^2}B_{b,g}\right\|_{t_k^{-1}}\left\|e^{-t_{k+1}tD_g^2}
\right\|_{t_{k+1}^{-1}}\\ \leq C^k t^{-{k\over 2}}\left(t_1\cdots
t_{k}\right)^{-{1\over 2}}\left({\rm Tr}\left[e^{-{tD_g^2\over
2}}\right]\right)^{t_1+\cdots+t_k}\left({\rm Tr}\left[e^{-{tD_g^2
}}\right]\right)^{t_{k+1}}\\
\leq C^k t^{-{k\over 2}}\left(t_1\cdots t_{k}\right)^{-{1\over 2}}
{\rm Tr}\left[e^{-{tD_g^2\over 2}}\right] .
\end{multline}
Thus for any $k\geq 1$, $t>0$,  one has
\begin{multline}\label{6.13}
 \left\|t^k\int_{\Delta_k} e^{-t_1tD_g^2}B_{b,g}
e^{-t_2tD_g^2}\cdots B_{b,g}e^{-t_{k+1}tD_g^2} dt_1\cdots
dt_k\right\|_1\\
\leq   \left(2C\sqrt{t}\right)^k{\rm Tr}\left[e^{-{tD_g^2\over
2}}\right]\int_{\Delta_k}d\sqrt{t_1}\cdots d\sqrt{t_k}\\ \leq
\left(2C\sqrt{t}\right)^k{\rm Tr}\left[e^{-{tD_g^2\over
2}}\right].
\end{multline}

From (\ref{6.4}) and (\ref{6.13}), one sees that at least for
$0<t\leq \min\{1,{1\over 8C^2}\}$, one has
\begin{align}\label{6.14}
 e^{-tD_b^2}=e^{-tD_g^2}+\sum_{k=1}^{+\infty}(-1)^kt^k\int_{\Delta_k}e^{-t_1tD_g^2}B_{b,g}
e^{-t_2tD_g^2}\cdots B_{b,g}e^{-t_{k+1}tD_g^2}dt_1\cdots dt_k.
\end{align}

From (\ref{6.6}), (\ref{6.13}) and (\ref{6.14}), one gets
(\ref{6.10}) easily.

The proof of Lemma \ref{t6.2} is completed. \ \ Q.E.D.

$\ $

From (\ref{6.8})-(\ref{6.10}) and proceed as in (\ref{6.12}) and
(\ref{6.13}), one deduces that when $t>0$ is small enough,
\begin{align}\label{6.15}
\left|t^{n+1}\int_{\Delta_{n+1}}{\rm
Tr}_s\left[Ne^{-t_1tD_g^2}B_{b,g} e^{-t_2tD_g^2}\cdots
B_{b,g}e^{-t_{n+2}tD_b^2}\right]dt_1\cdots dt_{n+1}\right|
 = O\left({  t^{1\over 2}}\right),
\end{align}
which completes the proof of Proposition \ref{t6.1}.\ \ Q.E.D.

$\ $

To compute the local index contribution to other terms in
(\ref{6.4}), we give the following formula for $B_{b,g}$.

\begin{thm}\label{t6.3} The following identity holds,
\begin{multline}\label{6.16}
 D_b^2=D_g^2+{1\over 2}\sum_{i,\, j=1}^n
  c\left(e_i\right)\widehat{c}\left(e_j\right)\left(\nabla_{e_i}^u\omega^F\left(e_j\right)\right)
\\  -{1\over 2}\sum_{i,\, j=1,\ i\neq j}^n
  c\left(e_i\right) {c}\left(e_j\right)\left(\nabla_{e_i}^u\omega^F\left(e_j\right)\right)
+{1\over 2}\sum_{i=1}^n
\left(\nabla_{e_i}^u\omega^F\left(e_i\right)\right)
  +\sum_{i=1}^n\omega^F\left(e_i\right)\nabla^u_{e_i}\\
  +{1\over
  4}\left(\widehat{c}\left(\omega^F\right)-
  {c}\left(\omega^F\right)\right)^2-{1\over
  4}\left[\widehat{c}\left(\omega^F\right)-
  {c}\left(\omega^F\right),\widehat{c}\left(\omega^F_g\right) \right].
\end{multline}
\end{thm}

{\it Proof}. From (\ref{4.8}) and (\ref{6.3}), one has
\begin{multline}\label{6.17}
  D_b^2-D_g^2=\left(d^F+d^{F*}_b\right)^2-\left(d^F+d^{F*}_g\right)^2\\
  =\left(d^F+d^{F*}_g+{1\over 2}\widehat{c}\left(\omega^F\right)
  -{1\over 2} {c}\left(\omega^F\right)
  \right)^2-\left(d^F+d^{F*}_g\right)^2\\
  ={1\over
  2}\left[d+d^{F*}_g,\widehat{c}\left(\omega^F\right)- {c}\left(\omega^F\right)
  \right] +{1\over
  4}\left(\widehat{c}\left(\omega^F\right)-
  {c}\left(\omega^F\right)\right)^2\\
  ={1\over
  2}\left[\sum_{i=1}^n
  c\left(e_i\right)\nabla_{e_i}^u-{1\over 2}\widehat{c}\left(\omega_g^F\right)
  ,\widehat{c}\left(\omega^F\right)- {c}\left(\omega^F\right)\right] +{1\over
  4}\left(\widehat{c}\left(\omega^F\right)-
  {c}\left(\omega^F\right)\right)^2.
\end{multline}

Now we compute (cf. \cite[4.33]{BZ1})
\begin{align}\label{6.18}
\left[\sum_{i=1}^n
  c\left(e_i\right)\nabla_{e_i}^u
  ,\widehat{c}\left(\omega^F\right) \right]=\sum_{i,\, j=1}^n
  c\left(e_i\right)\widehat{c}\left(e_j\right)\left(\nabla_{e_i}^u\omega^F\left(e_j\right)\right),
\end{align}
\begin{multline}\label{6.19}
\left[\sum_{i=1}^n
  c\left(e_i\right)\nabla_{e_i}^u
  , {c}\left(\omega^F\right) \right]=\sum_{i,\, j=1,\ i\neq j}^n
  c\left(e_i\right) {c}\left(e_j\right)\left(\nabla_{e_i}^u\omega^F\left(e_j\right)\right)
\\
-\sum_{i=1}^n \left(\nabla_{e_i}^u\omega^F\left(e_i\right)\right)
  -2\omega^F\left(e_i\right)\nabla^u_{e_i}.
\end{multline}

From (\ref{6.17})-(\ref{6.19}), we get (\ref{6.16}).\ \ Q.E.D.

$\ $

To compute the local index, let $a>0$ be the injectivity radius of
$(M,g^{TM})$. Take $x\in M$, let $e_1,\,\cdots,\,e_n$ be an
orthonormal basis of $T_xM$. We identify the open ball
$B^{T_xM}(0,a/2)$ with the open ball $B^M(x,a/2)$ in $M$ using
geodesic coordinates. Then $y\in T_xM,\ |y|\leq a/2$, represents
an element of $B^M(0,a/2)$. For $y\in T_xM$, $|y|\leq a/2$, we
identify $T_yM$, $F_y$ to $T_xM$, $F_x$ by parallel transport
along the geodesic $t\in[0,1]\rightarrow ty$ with respect to the
connections $\nabla^{TM}$, $\nabla^{F,u}$ respectively.

Let $\Gamma^{TM,x}$, $\Gamma^{F,u,x}$ be the connection forms for
$\nabla^{TM}$ $\nabla^{F,u}$ in the considered trivialization of
$TM$. By \cite[Proposition 4.7]{ABP}, one has
\begin{align}\label{6.20}
\Gamma^{TM,x}_y={1\over 2}R^{TM}_x(y,\cdot )+O\left(|y|^2\right),
\end{align}
$$ \Gamma^{F,u,x}_y=O(|y|).$$

Following \cite[(4.20)]{BZ1}, for any $t>0$, we   introduce the
Getzler rescaling
\begin{align}\label{6.21}
 {c}_t\left(e_i\right)= {e_i \over
t^{1/4}}\wedge-t^{1/4}i_{e_i} ,\ \ \ \
\widehat{c}_t\left(e_i\right)= {\widehat{e_i} \over
t^{1/4}}\wedge-t^{1/4}i_{\widehat{e_i}} ,\ \ \ y\rightarrow
\sqrt{t}y,
\end{align}
where we have written $e^*_i\wedge$ in \cite[(4.20)]{BZ1} as
$e_i\wedge$ for the sake of simplicity.

From (\ref{6.3}), (\ref{6.16}), one verifies easily that under the
Getzler rescaling  $G_t$ defined in (\ref{6.21}), one has that as
$t\rightarrow 0^+$,
\begin{multline}\label{6.22}
 G_t\left(tB_{b,g}\right)=\sqrt{t}\left( {1\over 2} \sum_{i,\, j=1}^n
   e_i \wedge\widehat{e_j}\left(\nabla_{e_i}^u\omega^F\left(e_j\right)\right)
   -{1\over 2}\sum_{i,\, j=1 }^n
   e_i \wedge   e_j
   \left(\nabla_{e_i}^u\omega^F\left(e_j\right)\right)\right.
  \\ +\left.\sum_{i=1}^n\omega^F\left(e_i\right){\partial
  \over\partial y^i}
  +{1\over
  4}\left(\widehat{ \omega^F}-
   \omega^F \right)^2-{1\over
  4}\left[\widehat{ \omega^F }-
   \omega^F ,\widehat{ \omega^F_g} \right]\right)+O(t).
\end{multline}

On the other hand, by \cite[(11.1)]{BZ1}, one has
\begin{align}\label{6.23}
  G_t(N)={1\over 2\sqrt{t}}\sum_{i=1}^ne_i\wedge \widehat{e_i}+O(1)
  ={1\over \sqrt{t}}L+O(1).
\end{align}

From (\ref{6.22}), (\ref{6.23}) and proceed as in \cite[Section
4]{BZ1}, \cite{G1} and \cite{G2}, one deduces that for any $ 1<
k\leq n$, $(t_1,\cdots,t_{k+1})\in\Delta_k$,
\begin{align}\label{6.24}
\lim_{t\rightarrow 0^+}t^k {\rm Tr}_s\left[Ne^{-t_1tD_g^2}B_{b,g}
e^{-t_2tD_g^2}\cdots B_{b,g}e^{-t_{k+1}tD_g^2}\right]=0,
\end{align}
while for $k=1$, $0\leq t_1\leq 1$,  one has
\begin{multline}\label{6.25}
\lim_{t\rightarrow 0^+}t {\rm Tr}_s\left[Ne^{-t_1tD_g^2}B_{b,g}
 e^{-(1-t_1)tD_g^2}\right]=\lim_{t\rightarrow 0^+}t {\rm Tr}_s\left[NB_{b,g}
 e^{- tD_g^2}\right]\\
 ={1\over 2}\int_M \int^B {\rm Tr}\left[  \sum_{i,\, j=1}^n
   e_i \wedge\widehat{e_j}\left(\nabla_{e_i}^u\omega^F\left(e_j\right)\right)+{1\over
  2}\left[
   \omega^F ,\widehat{ \omega^F_g}-\widehat{ \omega^F }
   \right]\right]L\exp\left(-{\dot{R}^{TM}\over 2}\right).
\end{multline}

Now it is clear that
\begin{align}\label{6.26}
 {\rm Tr}\left[
   \omega^F ,\widehat{ \omega^F_g}-\widehat{ \omega^F }
   \right] =0,
\end{align}
while by   \cite[(4.73)]{BZ1} and using the notation in
\cite[Section 4]{BZ1},
\begin{align}\label{6.27}
\sum_{i,\, j=1}^n
   e_i \wedge\widehat{e_j}{\rm
   Tr}\left[\left(\nabla_{e_i}^u\omega^F\left(e_j\right)\right)\right]
   = \nabla^{TM}\varphi  {\rm Tr}\left[\omega^F\right]
    ,
\end{align}
from which and from \cite[(3.10)]{BZ1} and \cite[(3.53)]{BZ1}, one
gets
\begin{multline}\label{6.28}
\int_M\int^B  \sum_{i,\, j=1}^n
   e_i
   \wedge\widehat{e_j}{\rm Tr}\left[\left(\nabla_{e_i}^u\omega^F\left(e_j\right)\right)\right]
  L\exp\left(-{\dot{R}^{TM}\over 2}\right)\\ =
\int_M\int^B  \nabla^{TM}\left(\left(\varphi  {\rm
Tr}\left[\omega^F\right]\right)
  L\exp\left(-{\dot{R}^{TM}\over 2}\right) \right)  =0.
    \end{multline}

    From (\ref{6.25}), (\ref{6.26}) and (\ref{6.28}), one gets for
    any $0\leq t_1\leq 1$ that
    \begin{align}\label{6.29}
\lim_{t\rightarrow 0^+}t {\rm Tr}_s\left[Ne^{-t_1tD_g^2}B_{b,g}
 e^{-(1-t_1)tD_g^2}\right]=0.
\end{align}

From (\ref{6.4}), (\ref{6.5}), (\ref{6.24}), (\ref{6.29}) and
\cite[Theorem 7.10]{BZ1}, one gets (\ref{3.16}).

The proof of Theorem \ref{t3.6} is completed. \ \ Q.E.D.

\begin{Rem}\label{6.4} The method developed in this
section, combined  with the method in \cite[Section 4]{BZ1}, can
be used to give an alternate  proof of Theorem \ref{t2.9}.
\end{Rem}

\section{Proof of Theorem \ref{t3.7}} \label{s7}
\setcounter{equation}{0}

We first restate Theorem \ref{t3.7} as follows.

\begin{thm}\label{t7.1}
There exist $0<d\leq 1$, $C>0$ such that for any $0<t\leq d$,
$0\leq T\leq {1\over t}$, then
\begin{multline}\label{7.1}
 \left|  {\rm Tr}_s\left[N\exp\left(-\left(tD_{b }+T\widehat{c}
 (\nabla f)\right)^2\right)\right]-
 {1\over t}\int_M\int^BL\exp\left(-B_{T^2}\right){\rm rk}(F)\right.\\
\left.  -{T\over
2}\int_M\theta\left(F,b^F\right)\int^B\widehat{df}\exp\left(-B_{T^2}\right)-
  {n\over 2}\chi(F)\right|\leq Ct.
\end{multline}
\end{thm}

Set, in view of (\ref{4.6}),
\begin{align}\label{7.2}
 \theta\left(F,g^F\right)={\rm Tr}\left[\omega^F_g\right]
 ={\rm Tr}\left[\left(g^F\right)^{-1}\nabla^Fg^F\right].
\end{align}

By \cite[Theorem A.1]{BZ2}, one has, under the same conditions as
in Theorem \ref{t7.1},
\begin{multline}\label{7.3}
 \left|  {\rm Tr}_s\left[N\exp\left(-\left(tD_{g }+T\widehat{c}
 (\nabla f)\right)^2\right)\right]-
 {1\over t}\int_M\int^BL\exp\left(-B_{T^2}\right){\rm rk}(F)\right.\\
\left.  -{T\over
2}\int_M\theta\left(F,g^F\right)\int^B\widehat{df}\exp\left(-B_{T^2}\right)-
  {n\over 2}\chi(F)\right|\leq C't
\end{multline}
for some constant $C'>0$.

Thus, in order to prove (\ref{7.1}), one need only to prove that
under the conditions of Theorem \ref{t7.1}, there exists constant
$C''>0$ such that
\begin{multline}\label{7.4}
 \left|  {\rm Tr}_s\left[N\exp\left(-\left(tD_{b }+T\widehat{c}
 (\nabla f)\right)^2\right)\right]- {\rm Tr}_s\left[N\exp\left(-\left(tD_{g}+T\widehat{c}
 (\nabla f)\right)^2\right)\right]\right.\\
\left.  -{T\over
2}\int_M\left(\theta\left(F,b^F\right)-\theta\left(F,g^F\right)\right)
\int^B\widehat{df}\exp\left(-B_{T^2}\right) \right|\leq C''t.
\end{multline}

Set for $t>0$, $T\geq 0$ that
\begin{align}\label{7.5}
  A_{b,t,T}=tD_b+T\widehat{c}(\nabla f),\ \ \ \ \ \  A_{g,t,T}=tD_g+T\widehat{c}(\nabla f)
\end{align}
and
\begin{align}\label{7.6}
  C_{t,T}=A_{b,t,T}^2-A_{g,t,T}^2 .
\end{align}
Then by (\ref{6.2}) and (\ref{6.3}) one has
\begin{multline}\label{7.7}
C_{t,T}=
  \left(tD_b+T\widehat{c}(\nabla f)\right)^2-\left(tD_g+T\widehat{c}(\nabla
  f)\right)^2
  =t^2B_{b,g}+tT\left[D_b-D_g,\widehat{c}(\nabla
  f)\right]\\
  =t^2B_{b,g}+{tT\over 2}\left[ \widehat{c}\left(\omega^F\right)-{c}\left(\omega^F\right)
  ,\widehat{c}(\nabla
  f)\right]=t^2B_{b,g}+{tT } \omega^F (\nabla f) .
\end{multline}

By (\ref{7.6}) and the Duhamel principle, one has
\begin{multline}\label{7.8}
 e^{-A_{b,t,T}^2}=e^{-A_{g,t,T}^2}\\
 +\sum_{k=1}^{n}(-1)^k\int_{\Delta_k}e^{-t_1A_{g,t,T}^2}C_{t,T}
e^{-t_2A_{g,t,T}^2}\cdots C_{t,T}e^{-t_{k+1}A_{g,t,T}^2}dt_1\cdots dt_k\\
+ (-1)^{n+1} \int_{\Delta_{n+1}}e^{-t_1A_{g,t,T}^2}C_{t,T}
e^{-t_2A_{g,t,T}^2}\cdots C_{t,T}e^{-t_{n+2}A_{b,t,T}^2}dt_1\cdots
dt_{n+1}.
\end{multline}

\begin{lemma}\label{t7.2} There exists $C_0>0$ such that for any
$T\geq 0$, $s\in \Omega^*(M,F)$, one has
\begin{align}\label{7.9}
 \left\|B_{b,g}s\right\|_0^2\leq C_0\left(\|s\|_0^2+
 \left\|\left(D_g+T\widehat{c}(\nabla f)\right)s\right\|_0^2\right).
\end{align}
\end{lemma}

{\it Proof}. Since both $b^F$ and $g^F$ by assumption are flat
near the set $B$ of critical points of the Morse function $f$, by
(\ref{6.3}) and (\ref{6.16}) we find that there exists $\delta>0$
such that
\begin{align}\label{7.10}
  B_{b,g}=0
\end{align}
on $\cup_{x\in B}B_x^M(2\delta)$, where for each $x\in B$,
$B_x^M(2\delta)\subset M$ is the ball of radius $2\delta$ centered
at $x$.

Let $\psi\geq 0$ be a function on $M$ such that ${\rm
supp}(\psi)\subset M\setminus\cup_{x\in B}B_x^M(\delta)$ while
$\psi\equiv 1$ on $M\setminus\cup_{x\in B}B_x^M({3\over
2}\delta)$. Then by (\ref{7.10}) and the standard elliptic
estimate, there exists $C_1>0$ such that for any $s\in
\Omega^*(M,F)$,
\begin{align}\label{7.11}
 \left\|B_{b,g}s\right\|_0^2 =
 \left\|B_{b,g}( \psi s)\right\|_0^2\leq C_1
 \left(\| \psi s\|_0^2+\left\|D_g( \psi s)\right\|^2_0\right) .
\end{align}

Also, by (\ref{4.18}) and (\ref{4.19}) it is clear that there
exists $C_2>0$ such that for any $T\geq 0$ and  $x\in
M\setminus\cup_{x\in B}B_x^M(\delta)$,
\begin{align}\label{7.12}
   T\left[D_g,\widehat{c}(\nabla f)\right]+T^2|\nabla f|^2 \geq -C_2 .
\end{align}

From (\ref{7.11}) and (\ref{7.12}), one deduces  that there exists
$C_3>0$ such that  for any $T\geq 0$ and any $s\in \Omega^*(M,F)$,
one has
\begin{multline}\label{7.13}
\left\|D_g( \psi s)\right\|^2_0\leq C_2\|\psi s\|_0^2 +
\left\langle \left(D_g+T\widehat{c}(\nabla f)\right)^2 (\psi
s),\psi s\right\rangle_g\\
= C_2\|\psi s\|_0^2 + \left\|\left(D_g+T\widehat{c}(\nabla
f)\right) (\psi s) \right\|_0^2\\
\leq C_3\left(\| s\|_0^2 + \left\|\left(D_g+T\widehat{c}(\nabla
f)\right)   s \right\|_0^2\right).
\end{multline}

From (\ref{7.11}) and (\ref{7.13}), one gets (\ref{7.9}). \ \
Q.E.D.

$\ $

By (\ref{7.5}), Lemma \ref{t7.2} and proceed  as in \cite[Lemma
2.8]{CH} and \cite[Lemma 1]{F}, one finds that there exists
$C_4>0$ such that for any $t>0,\ u>0$ verifying $ut^2\leq 1$ and $
T\geq 0$,
\begin{align}\label{7.14}
\left\| \exp\left(-uA_{g,t,T}^2\right)B_{b,g}\right\|_{u^{-1}}
\leq
 C_4 u^{-{1\over 2}}t^{-1}\left({\rm Tr}
 \left[\exp\left(-{A_{g,t,T}^2\over 2}\right)\right]\right)^u.
\end{align}

Similarly, as
\begin{align}\label{7.15}
  \omega^F =0
\end{align}
on $\cup_{x\in B}B_x^M(2\delta)$, one deduces that there exists
$C_5>0$ such that for any $u>0$, $t>0$, $T\geq 0$,
\begin{align}\label{7.16}
\left\| \exp\left(-uA_{g,t,T}^2\right) T  \omega^F (\nabla
  f) \right\|_{u^{-1}}  \leq
 C_5 u^{-{1\over 2}} \left({\rm Tr}
 \left[\exp\left(-{A_{g,t,T}^2\over 2}\right)\right]\right)^u.
\end{align}

From (\ref{6.8}), (\ref{6.9}), (\ref{7.7}), (\ref{7.14}),
(\ref{7.16}) and proceed as in (\ref{6.12}), one sees that for any
$k\geq 1$ and $t>0$, $t_i>0$ for $1\leq i\leq k+1$ with
$\sum_{i=1}^{k+1}t_i=1$, one has
\begin{multline}\label{7.17}
  \left\|e^{-t_1A_{g,t,T}^2}C_{t,T}
e^{-t_2A_{g,t,T}^2}\cdots C_{t,T}e^{-t_{k+1}A_{g,t,T}^2}\right\|_1
\\ \leq \left(C_4+C_5\right)^k t^{{k }} \left(t_1\cdots t_{t_k}\right)^{-{1\over
2}} {\rm Tr}\left[e^{-{A_{g,t,T}^2\over 2}}\right] .
\end{multline}

From (\ref{7.17}) and proceed as in (\ref{6.13}), one has for any
$k\geq 1$ and $t>0$ that
\begin{multline}\label{7.18}
  \left\|\int_{\Delta_k}  e^{-t_1A_{g,t,T}^2}C_{t,T}
e^{-t_2A_{g,t,T}^2}\cdots C_{t,T}e^{-t_{k+1}A_{g,t,T}^2}
dt_1\cdots dt_k\right\|_1
\\ \leq  \left(2\left(C_4+C_5\right)t\right)^k{\rm Tr}\left[e^{-{A_{g,t,T}^2\over 2}}\right]   .
\end{multline}

From (\ref{7.8}) and (\ref{7.18}), one gets that at least for
$0<t\leq \min\{1, {1\over 4(C_4+C_5)}\}$ and $T\geq 0$ with
$tT\leq 1$, one has
\begin{multline}\label{7.19}
 e^{-A_{b,t,T}^2}=e^{-A_{g,t,T}^2}\\
 +\sum_{k=1}^{+\infty}(-1)^k\int_{\Delta_k}e^{-t_1A_{g,t,T}^2}C_{t,T}
e^{-t_2A_{g,t,T}^2}\cdots C_{t,T}e^{-t_{k+1}A_{g,t,T}^2}dt_1\cdots
dt_k .
\end{multline}

From (6.6), (\ref{7.18}), (\ref{7.19}) and \cite[(12.34)]{BZ1},
one finds that for any $0<t\leq \min\{1, {1\over 4(C_4+C_5)}\}$
and $T\geq 0$ with $tT\leq 1$, one has that for any $r>0$, there
exists $C_6>0$ such that
\begin{align}\label{7.20}
\left\|\exp\left(-A_{b,t,T}^2\right)\right\|_r\leq {C_6\over t^{n
}}.
\end{align}

From (\ref{7.14}), (\ref{7.16}) and (\ref{7.20}), one can proceed
as in (\ref{6.12}) and (\ref{6.15}) to see that there exists
$C_7>0$ such that for any $t>0$ small enough and $T\in[0,{1\over
t}]$, one has
\begin{align}\label{7.21}
  \left|\int_{\Delta_{n+1}}{\rm Tr}_s\left[ Ne^{-t_1A_{g,t,T}^2}C_{t,T}
e^{-t_2A_{g,t,T}^2}\cdots
C_{t,T}e^{-t_{n+2}A_{b,t,T}^2}\right]dt_1\cdots dt_{n+1}\right|
\leq  {C_7 t }   .
\end{align}

Now for any $x\in M$, we introduce the coordinates and
identification around $x$ as in Section \ref{s6}, and use the
Getzler rescaling introduced in (\ref{6.21}), with $t$ there
replaced by $t^2$ here. By using    (\ref{7.7}), one has
\begin{align}\label{7.22}
  G_{t^2} \left(C_{t,T}\right)
  = G_{t^2}\left(t^2B_{b,g}\right)+ tT\omega^F(\nabla f).
\end{align}

From (\ref{6.21})-(\ref{6.29}), (\ref{7.22}) and proceed as in
\cite[Section 13]{BZ1}, one deduces that there exists $C_8>0$,
$0<d\leq 1$ such that for any $1<k\leq n$, $0<t\leq d$, $T\geq 0$
with $tT\leq 1$,
\begin{align}\label{7.23}
\left|\int_{\Delta_k}{\rm Tr}_s\left[Ne^{-t_1A_{g,t,T}^2}C_{t,T}
e^{-t_2A_{g,t,T}^2}\cdots
C_{t,T}e^{-t_{k+1}A_{g,t,T}^2}\right]dt_1\cdots dt_k\right|\leq
C_8t,
\end{align}
while for $k=1$ one has for any $0<t\leq d$, $T\geq 0$ with
$tT\leq 1$ and $0\leq t_1\leq 1$,
\begin{multline}\label{7.24}
\left| {\rm Tr}_s\left[Ne^{-t_1A_{g,t,T}^2}C_{t,T}
e^{-\left(1-t_1\right)A_{g,t,T}^2} \right] - T\int_M\int^B{\rm
Tr}\left[\omega^F(\nabla
f)\right]L\exp\left(-B_{T^2}\right)\right|\\ \leq C_8t.
\end{multline}

Now from \cite[(3.9)]{BZ1}, \cite[(3.52)]{BZ1},
\cite[(3.53)]{BZ1}, (\ref{2.26}), (\ref{4.2}), (\ref{6.1}) and
(\ref{7.2}), one deduces
  that
\begin{multline}\label{7.25}
 \int_M\int^B{\rm
Tr}\left[\omega^F(\nabla f)\right]L\exp\left(-B_{T^2}\right)=
\int_M\int^Bi_{ \nabla f}\left({\rm
Tr}\left[\omega^F\right]\right)L\exp\left(-B_{T^2}\right)\\
=\int_M\int^B {\rm Tr}\left[\omega^F\right] i_{ \nabla
f}(L)\exp\left(-B_{T^2}\right)+\int_M\int^B {\rm
Tr}\left[\omega^F\right] L\,i_{ \nabla
f}\left(\exp\left(-B_{T^2}\right)\right)\\
={1\over 2}\int_M\int^B {\rm Tr}\left[\omega^F\right]  \widehat{
\nabla f} \exp\left(-B_{T^2}\right)-{1\over 2}\int_M\int^B {\rm
Tr}\left[\omega^F\right]
L\,\nabla^{TM}\left(\exp\left(-B_{T^2}\right)\right)\\
={1\over 2}\int_M {\rm Tr}\left[\omega^F\right] \int^B \widehat{
\nabla f} \exp\left(-B_{T^2}\right)-{1\over 2}\int_M{\rm
Tr}\left[\omega^F\right]\int^B\nabla^{TM}\left( L
\left(\exp\left(-B_{T^2}\right)\right)\right)\\
={1\over 2}\int_M
\left(\theta\left(F,g^F\right)-\theta\left(F,b^F\right)\right)
\int^B \widehat{ \nabla f} \exp\left(-B_{T^2}\right).
\end{multline}

From (\ref{7.8}), (\ref{7.21}) and (\ref{7.23})-(\ref{7.25}), one
gets (\ref{7.4}), which completes the proof of Theorem \ref{t7.1}.
\ \ Q.E.D.

\section{Proof of Theorem \ref{t3.8}} \label{s8}
\setcounter{equation}{0}

In view of (\ref{3.18}) and \cite[Theorem A.2]{BZ2}, in order to
prove Theorem \ref{t3.8}, we need only to prove that for any
$T>0$,
\begin{align}\label{8.1}
\lim_{t\rightarrow 0^+}\left({\rm
Tr}_s\left[N\exp\left({-A_{b,t,{T\over t}}^2}\right)\right] -{\rm
Tr}_s\left[N\exp\left({-A_{g,t,{T\over
t}}^2}\right)\right]\right)=0.
\end{align}

First of all, by (\ref{7.18}), there exists $0<C_0\leq 1$ such
that when $0<t\leq C_0$, one has
\begin{multline}\label{8.2}
  \left\|\int_{\Delta_k}  e^{-t_1A_{g,t,{T\over t}}^2}C_{t,{T\over t}}
e^{-t_2A_{g,t,{T\over t}}^2}\cdots C_{t,{T\over
t}}e^{-t_{k+1}A_{g,t,{T\over t}}^2} dt_1\cdots dt_k\right\|_1
\\ \leq  \left( {t\over 2C_0}\right)^k{\rm Tr}\left[e^{-{A_{g,t,{T\over t}}^2\over 2}}\right]   .
\end{multline}
Thus we have the  absolute convergent expansion formula
\begin{multline}\label{8.3}
 e^{-A_{b,t,{T\over t}}^2}-e^{-A_{g,t,{T\over t}}^2}\\
 =\sum_{k=1}^{+\infty}(-1)^k\int_{\Delta_k}e^{-t_1A_{g,t,{T\over t}}^2}C_{t,{T\over t}}
e^{-t_2A_{g,t,{T\over t}}^2}\cdots C_{t,{T\over
t}}e^{-t_{k+1}A_{g,t,{T\over t}}^2}dt_1\cdots dt_k .
\end{multline}

Since $T>0$ is fixed, by \cite[(12.34) and (15.22)]{BZ1}, there
exists $C_1>0$ such that for $0<t\leq C_0$,
\begin{align}\label{8.4}
{\rm Tr}\left[e^{-{A_{g,t,{T\over t}}^2\over 2}}\right]\leq
{C_1\over t^n} .
\end{align}

From (\ref{8.2}) and (\ref{8.4}), one sees that
\begin{align}\label{8.5}
\sum_{k=n}^{+\infty}(-1)^k\int_{\Delta_k}e^{-t_1A_{g,t,{T\over
t}}^2}C_{t,{T\over t}} e^{-t_2A_{g,t,{T\over t}}^2}\cdots
C_{t,{T\over t}}e^{-t_{k+1}A_{g,t,{T\over t}}^2}dt_1\cdots dt_k
\end{align}
is uniformly absolute convergent for $0<t\leq C_0$.

Let $\psi\geq 0$ be the function on $M$ defined in Section
\ref{s7}. Then by definition  one has
\begin{align}\label{8.6}
C_{t,{T\over t}}=\psi C_{t,{T\over t}}=C_{t,{T\over t}}\psi =\psi
C_{t,{T\over t}}\psi .
\end{align}

From (\ref{8.3}) one sees that for each $k\geq 1$ and any $T>0$,
$0<t\leq C_0$ and $(t_1,\cdots, t_{k+1})\in \Delta_k$, since
$\sum_{i=1}^{k+1}t_i=1$,  there is $j\in [1,k+1]$ such that
$t_j\geq {1\over k+1}$. We here deal with the case where $j=k+1$,
the other cases can be dealt with similarly.

From (\ref{6.8}), (\ref{6.9}), (\ref{7.7}), (\ref{7.14}),
(\ref{7.16}), (\ref{8.6}) and proceed as in (\ref{6.12}),  one has
that for any $(t_1,\cdots, t_{k+1})\in\Delta_k\setminus
\{t_1\cdots t_{k+1}=0\}$,
\begin{multline}\label{8.7}
\left|{\rm Tr}_s\left[N  e^{-t_1A_{g,t,{T\over t}}^2}C_{t,{T\over
t}} e^{-t_2A_{g,t,{T\over t}}^2}\cdots C_{t,{T\over
t}}e^{-t_{k+1}A_{g,t,{T\over t}}^2}\right]\right| \\ =\left|{\rm
Tr}_s\left[N e^{-t_1A_{g,t,{T\over t}}^2}C_{t,{T\over t}}
e^{-t_2A_{g,t,{T\over t}}^2}\cdots C_{t,{T\over t}}\psi
e^{-t_{k+1}A_{g,t,{T\over t}}^2}\right]\right|
 \\ \leq C_2 \left\|
e^{- t_1 A_{g,t,{T\over t}}^2} C_{t,{T\over t}}\right\|_{ t_1
^{-1} }   \left\|e^{-{t_2 }A_{g,t,{T\over t}}^2} C_{t,{T\over
t}}\right\|_{t_2^{-1}} \\ \cdots \left\| e^{-t_kA_{g,t,{T\over
t}}^2} C_{t,{T\over t}}\right\|_{t_k^{-1}}  \left\|\psi
e^{-t_{k+1}A_{g,t,{T\over t}}^2}  \right\|_{t_{k+1}^{-1}}
\\ \leq  C_3t^k  \left(t_1\cdots t_{k}\right)^{-{1\over
2}} {\rm Tr}\left[e^{-{A_{g,t,{T\over t}}^2\over 2}}\right]
\left\| \psi e^{-{ t_{k+1}\over 2}A_{g,t,{T\over t}}^2} \right\|
\end{multline}
for some positive constants $C_2>0,\ C_3>0$.

From (\ref{8.4}),   (\ref{8.7})  and the assumption that
$t_{k+1}\geq {1\over k+1}$, one gets
\begin{multline}\label{8.8}
\left|\int_{\Delta_k}{\rm Tr}_s\left[N  e^{-t_1A_{g,t,{T\over
t}}^2}C_{t,{T\over t}} e^{-t_2A_{g,t,{T\over t}}^2}\cdots
C_{t,{T\over t}}e^{-t_{k+1}A_{g,t,{T\over t}}^2}\right]dt_1\cdots dt_k\right| \\
 \leq C_4  t^{k-n}    \left\| \psi
e^{-{ 1\over 2(k+1)}A_{g,t,{T\over t}}^2} \right\|
\end{multline}
for some constant $C_4>0$.

By   (\ref{8.8}) one need to estimate
\begin{multline}\label{8.9}
   \left\| \psi
e^{-{ 1\over 2(k+1)}A_{g,t,{T\over t}}^2} \right\|=\left\|\psi
e^{-{ 1\over 2(k+1)}A_{g,t,{T\over t}}^2}\left(\psi e^{-{ 1\over
2(k+1)}A_{g,t,{T\over t}}^2}\right)^*\right\|^{1\over 2} \\
=\left\|\psi e^{-{ 1\over k+1}A_{g,t,{T\over t}}^2}\psi
\right\|^{1\over 2}\\ \leq\sqrt{\int_M{\rm
Tr}\left[\psi(x)S_{{t\over \sqrt{  k+1 }},{1\over \sqrt{ k+1
}}{T\over t}}(x,x)\psi(x)\right]d{\rm vol}_x},
\end{multline}
where as in \cite[Section 14]{BZ1}, $S_{t,{T\over t}}(x,y)$ for
$x,\ y\in M$ denotes the kernel of $\exp(-A_{g,t,{T\over t}}^2)$
with respect to the Riemannian volume $d{\rm vol}_{g^{TM}}$.

Now since ${\rm Supp}(\psi)\subset M\setminus \cup_{x\in B}
B_x(\delta)$, by \cite[Proposition 14.1]{BZ1}, one sees that there
exist $C_5,\ C_6>0$ such that
\begin{align}\label{8.10}
\int_M{\rm Tr}\left[\psi(x)S_{{t\over \sqrt{  k+1 }},{1\over
\sqrt{ k+1 }}{T\over t}}(x,x)\psi(x)\right]d{\rm vol}_x\leq
C_5\exp\left(-{C_6\over t^2}\right).
\end{align}

From (\ref{8.3}), (\ref{8.5}), (\ref{8.8})-(\ref{8.10}) and the
dominate convergence, we get (\ref{8.1}), which completes the
proof of Theorem \ref{3.8}. \ \ Q.E.D.

\section{Proof of Theorem \ref{t3.9}} \label{s9}
\setcounter{equation}{0}

In view of (\ref{3.19}) and \cite[Theorem A.3]{BZ2}, in order to
prove Theorem \ref{t3.9}, we need only to prove that there exist
$c>0$, $C>0$, $0<C_0\leq 1$ such that for any  $0<t\leq C_0$, $
T\geq 1$,
\begin{align}\label{9.1}
 \left| {\rm
Tr}_s\left[N\exp\left({-A_{b,t,{T\over t}}^2}\right)\right] -{\rm
Tr}_s\left[N\exp\left({-A_{g,t,{T\over
t}}^2}\right)\right]\right|\leq c\exp(-CT).
\end{align}

 First of all, one can choose $C_0>0$ small enough so that  for any $0<t\leq C_0$, $T>0$, by (\ref{8.3}),
  we have the  absolute convergent expansion formula
\begin{multline}\label{9.2}
 e^{-A_{b,t,{T\over t}}^2}-e^{-A_{g,t,{T\over t}}^2}\\
 =\sum_{k=1}^{+\infty}(-1)^k\int_{\Delta_k}e^{-t_1A_{g,t,{T\over t}}^2}C_{t,{T\over t}}
e^{-t_2A_{g,t,{T\over t}}^2}\cdots C_{t,{T\over
t}}e^{-t_{k+1}A_{g,t,{T\over t}}^2}dt_1\cdots dt_k ,
\end{multline}
from which one has
\begin{multline}\label{9.3}
 {\rm
Tr}_s\left[N\exp\left({-A_{b,t,{T\over t}}^2}\right)\right] -{\rm
Tr}_s\left[N\exp\left({-A_{g,t,{T\over
t}}^2}\right)\right]\\
 =\sum_{k=1}^{+\infty}(-1)^k\int_{\Delta_k}{\rm Tr}_s\left[N
 e^{-t_1A_{g,t,{T\over t}}^2}C_{t,{T\over t}}
e^{-t_2A_{g,t,{T\over t}}^2}\cdots C_{t,{T\over
t}}e^{-t_{k+1}A_{g,t,{T\over t}}^2}\right]dt_1\cdots dt_k .
\end{multline}
Thus, in order to prove (\ref{9.1}), we need only to prove
\begin{multline}\label{9.4}
  \sum_{k=1}^{+\infty}\left|\int_{\Delta_k}{\rm Tr}_s\left[N
 e^{-t_1A_{g,t,{T\over t}}^2}C_{t,{T\over t}}
e^{-t_2A_{g,t,{T\over t}}^2}\cdots C_{t,{T\over
t}}e^{-t_{k+1}A_{g,t,{T\over t}}^2}\right]dt_1\cdots dt_k\right|\\
=\sum_{k=1}^{+\infty}\left|\int_{\Delta_k}{\rm Tr}_s\left[N
 e^{-\left(t_1+t_{k+1}\right)A_{g,t,{T\over t}}^2}C_{t,{T\over t}}
e^{-t_2A_{g,t,{T\over t}}^2}\cdots C_{t,{T\over
t}}\right]dt_1\cdots dt_k\right|\\ \leq c\exp(-CT) .
\end{multline}

Let $\psi\geq 0$ be the function on $M$ defined in Section
\ref{s7}. By (\ref{8.6}), we have for any $t>0$, $T\geq 1$,
$(t_1,\cdots,t_{k+1})\in \Delta_k\setminus\{t_1\cdots
t_{k+1}=0\}$,
\begin{multline}\label{9.5}
    {\rm Tr}_s\left[N
 e^{-\left(t_1+t_{k+1}\right)A_{g,t,{T\over t}}^2}C_{t,{T\over t}}
e^{-t_2A_{g,t,{T\over t}}^2}\cdots C_{t,{T\over
t}}\right] \\
={\rm Tr}_s\left[N
 \psi e^{-\left(t_1+t_{k+1}\right)A_{g,t,{T\over t}}^2}  C_{t,{T\over t}}
\psi e^{-t_2A_{g,t,{T\over t}}^2}  C_{t,{T\over t}}\cdots\psi
e^{-t_kA_{g,t,{T\over t}}^2}  C_{t,{T\over t}}\right] .
\end{multline}

We first state a refinement  of the estimates (\ref{6.11}),
(\ref{7.14}) and (\ref{7.16}).

\begin{lemma}\label{t9.1}
There exists $C_1>0$ such that for any  $0<u\leq 1$,     $0<t\leq
1$, $T\geq 1$, one has
\begin{align}\label{9.6}
\left\| \psi e^{-uA_{g,t,{T\over t}}^2}  C_{t,{T\over
t}}\right\|_{u^{-1}}\leq C_1u^{-{1\over 2}}t \left({\rm
Tr}\left[e^{-{1\over 2}A_{g,t,{T\over
t}}^2}\right]\right)^u\left\|\psi e^{-{u\over 4}A_{g,t,{T\over
t}}^2} \right\|.
\end{align}
\end{lemma}

{\it Proof}. From (\ref{7.9}), (\ref{7.15}) and (\ref{7.16}), one
has that there exists constant $C_2>0$
\begin{align}\label{9.7}
   C_{t,{T\over t}}^*  C_{t,{T\over t}} \leq C_2 t^2 \left(1+A_{g,t,{T\over
t}}^2\right) .
\end{align}

From (\ref{6.8}) and (\ref{9.7}), one gets
\begin{multline}\label{9.8}
\left\| \psi e^{-uA_{g,t,{T\over t}}^2}  C_{t,{T\over
t}}\right\|_{u^{-1}}\leq \left\| \psi e^{-{3u\over
4}A_{g,t,{T\over t}}^2}\right\|_{u^{-1}} \left\|e^{-{u\over
4}A_{g,t,{T\over t}}^2}C_{t,{T\over t}}\right\|   \\
\leq C_3u^{-{1\over 2}}t \left\|\psi e^{-{u\over 4}A_{g,t,{T\over
t}}^2}\right\|\left\|  e^{-{u\over 2}A_{g,t,{T\over
t}}^2}\right\|_{u^{-1}}\left\|e^{-{u\over 4}A_{g,t,{T\over t}}^2}
 \left(1+{u\over 4}A_{g,t,{T\over t}}^2\right)^{1\over 2}\right\|\\
\leq C_4 u^{-{1\over 2}}t \left({\rm Tr}\left[e^{-{1\over
2}A_{g,t,{T\over t}}^2}\right]\right)^u\left\|\psi e^{-{u\over
4}A_{g,t,{T\over t}}^2} \right\|
\end{multline}
for some positive constants $C_3>0$, $C_4>0$.

The proof of Lemma \ref{t9.1} is completed. \ \ Q.E.D.

$\ $

\begin{lemma}\label{t9.2}
There exist $0<c_1\leq 1$, $C_5>0$, $C_6>0$  such that for any $0<
u\leq 1$, $0<t\leq c_1$, $T\geq 1$, one has
\begin{align}\label{9.9}
\left\|\psi e^{-{u }A_{g,t,{T\over t}}^2} \right\|\leq
C_5\exp\left(-C_6uT\right) .
\end{align}
\end{lemma}

{\it Proof}. From (\ref{4.13}) and (\ref{7.5}), one has
\begin{align}\label{9.10}
 A_{g,t,{T\over t}}^2=t^2\left(D_g+{T\over t^2}\widehat{c}(\nabla
 f)\right)^2
 =t^2\widetilde{D}_{g,{T\over t^2}} ^2.
\end{align}

Since $T\geq 1$, it is known (cf. \cite{S1}) that there exists
$0<c_1\leq 1$, $c_2>0$, $c_3>0$  such that for any $0<t\leq c_1$,
the spectrum of $\widetilde{D}_{g,{T\over t^2}}^2 $ splits into
two parts:
\begin{align}\label{9.11}
{\rm Spec}\left(\widetilde{D}_{g,{T\over t^2}}^2\right)\subset
\left[0,
 \exp\left(-{c_2T\over t^2}\right)\right]\bigcup \left[{c_3T\over t^2},+\infty\right) .
\end{align}

For $0<t\leq c_1$ and $T\geq 1$ let $Q^{[0,1]}_{T\over t^2}$
denote the orthogonal projection from $L^2(\Omega^*(M,F)$ to the
direct sum of the eigenspaces of $\widetilde{D}_{g,{T\over
t^2}}^2$ corresponding to the eigenvalues lying in $[0,1]$. Let
$Q^{[1,+\infty)}_{T\over t^2}={\rm Id}- Q^{[0,1]}_{T\over t^2}$.
Then it is known that (cf. \cite[(7.20)]{BZ1}) ${\rm
Im}(Q^{[0,1]}_{T\over t^2})$ is a finite dimensional space.

Now we write
\begin{multline}\label{9.12}
\left\|\psi e^{-{u }A_{g,t,{T\over t}}^2} \right\| =\left\|\psi
e^{-{u }A_{g,t,{T\over t}}^2} \left(Q^{[0,1]}_{T\over
t^2}+Q^{[1,+\infty)}_{T\over t^2}\right) \right\|\\
\leq \left\|\psi e^{-{u }A_{g,t,{T\over t}}^2}
 Q^{[0,1]}_{T\over t^2}\right\|+\left\|\psi e^{-{u }A_{g,t,{T\over t}}^2}Q^{[1,+\infty)}_{T\over
 t^2}
\right\|.
\end{multline}

From (\ref{9.10}) and (\ref{9.11}), one sees that
\begin{align}\label{9.13}
\left\|\psi e^{-{u }A_{g,t,{T\over t}}^2}Q^{[1,+\infty)}_{T\over
 t^2}
\right\|\leq \left\|  e^{-{u }A_{g,t,{T\over
t}}^2}Q^{[1,+\infty)}_{T\over
 t^2}
\right\|\leq \exp\left(-c_3uT\right) .
\end{align}

From (\ref{9.11}) one has
\begin{multline}\label{9.14}
\left\|\psi e^{-{u }A_{g,t,{T\over t}}^2}
 Q^{[0,1]}_{T\over t^2}\right\|\leq  \left\|\psi \left( e^{-{u }A_{g,t,{T\over
 t}}^2}-{\rm Id}\right)
 Q^{[0,1]}_{T\over t^2}\right\|+\left\|\psi
 Q^{[0,1]}_{T\over t^2}\right\|\\
 \leq \left\|\psi
 Q^{[0,1]}_{T\over t^2}\right\|+ C_7\exp\left(-{c_4T\over
 t^2}\right)
\end{multline}
for some positive constants $c_4>0$, $C_7>0$.

For any $T>0$, let $J_T$ be the map defined in (\ref{4.31}) where
we assume without loss of generality that the radius $4a$ there
verifies $4a\leq \delta$. Then one has
\begin{align}\label{9.15}
\psi J_{T\over t^2}=0.
\end{align}

By (\ref{9.15}) and \cite[Theorem 8.8]{BZ1} and \cite[Theorem
6.7]{BZ2}, an analogue of which has been proved in Theorem
\ref{t4.3}, one sees that there exist $C_8>0$, $c_5>0$ such that
\begin{align}\label{9.16}
\left\|\psi Q^{[0,1]}_{T\over t^2} J_{T\over t^2}\right\|\leq
C_8\exp\left(-{c_5T\over
 t^2}\right) .
\end{align}

From (\ref{9.16}) one deduces easily that there exist $C_9>0$,
$c_6>0$ such that
\begin{align}\label{9.17}
\left\|\psi Q^{[0,1]}_{T\over t^2} \right\|\leq
C_9\exp\left(-{c_6T\over
 t^2}\right) .
\end{align}

From (\ref{9.12})-(\ref{9.14}) and (\ref{9.17}), one gets
(\ref{9.9}).

The proof of Lemma \ref{t9.2} is completed.\ \ Q.E.D.

$\ $

From Lemmas \ref{t9.1} and \ref{t9.2}, one deduces that for any
$0<t\leq \min\{C_0,\ c_1\}$, $T\geq 1$ and
$(t_1,\cdots,t_{k+1})\in \Delta_k\setminus\{t_1\cdots t_{k+1}=0\}$
that
\begin{multline}\label{9.18}
\left\|\psi e^{-\left(t_1+t_{k+1}\right)A_{g,t,{T\over t}}^2}
C_{t,{T\over t}} \psi e^{-t_2A_{g,t,{T\over t}}^2}  C_{t,{T\over
t}}\cdots\psi e^{-t_kA_{g,t,{T\over t}}^2}  C_{t,{T\over
t}}\right\|_1\\
\leq \left\| \psi e^{-\left(t_1+t_{k+1}\right)A_{g,t,{T\over
t}}^2} C_{t,{T\over t}}
\right\|_{\left(t_1+t_{k+1}\right)^{-1}}\cdots \left\| \psi e^{-
t_{k } A_{g,t,{T\over
t}}^2} C_{t,{T\over t}} \right\|_{  t_{k }^{-1}}\\
\leq \left(C_1C_5t\right)^k
\left(\left(t_1+t_{k+1}\right)t_2\cdots t_k\right)^{-{1\over
2}}{\rm Tr}\left[e^{-{1\over 2}A_{g,t,{T\over
t}}^2}\right]\exp\left(-{C_6T\over 4}\right).
\end{multline}

By \cite[(15.22)]{BZ1}, one sees that there exists $C_{10}>0$ such
that for any $0<t\leq \min\{C_0,\ c_1\}$, $T\geq 1$,
\begin{align}\label{9.19}
{\rm Tr}\left[e^{-{1\over 2}A_{g,t,{T\over t}}^2}\right]\leq
C_{10}{T^{n\over 2}\over t^n}.
\end{align}

From (\ref{9.5}), (\ref{9.18}) and (\ref{9.19}), one sees that
there exists $C_{11}>0$ such that for any $k\geq 1$,
\begin{multline}\label{9.20}
 \left|\int_{\Delta_k}{\rm
Tr}_s\left[N
 e^{-t_1A_{g,t,{T\over t}}^2}C_{t,{T\over t}}
e^{-t_2A_{g,t,{T\over t}}^2}\cdots C_{t,{T\over
t}}e^{-t_{k+1}A_{g,t,{T\over t}}^2}\right]dt_1\cdots dt_k\right|\\
\leq C_{11} \left(2\,C_1C_5t\right)^k{T^{n\over 2}\over t^n}
  \exp\left(-{C_6T\over 4}\right),
\end{multline}
from which one sees that there exist $0<c_7\leq \min\{C_0,\
c_1\}$, $C_{12}>0$, $C_{13}>0$ such that for any $0<t\leq c_7$ and
$T\geq 1$, one has
\begin{multline}\label{9.21}
 \left|\sum_{k=n}^{+\infty}\int_{\Delta_k}{\rm
Tr}_s\left[N
 e^{-t_1A_{g,t,{T\over t}}^2}C_{t,{T\over t}}
e^{-t_2A_{g,t,{T\over t}}^2}\cdots C_{t,{T\over
t}}e^{-t_{k+1}A_{g,t,{T\over t}}^2}\right]dt_1\cdots dt_k\right|\\
\leq C_{12}
  \exp\left(-C_{13}T\right).
\end{multline}

On the other hand, for any $1\leq k<n$, by proceeding as in
(\ref{8.8}), one has that for any $0<t\leq c_7$, $T\geq 1$,
\begin{multline}\label{9.22}
\left|\int_{\Delta_k}{\rm Tr}_s\left[N  e^{-t_1A_{g,t,{T\over
t}}^2}C_{t,{T\over t}} e^{-t_2A_{g,t,{T\over t}}^2}\cdots
C_{t,{T\over t}}e^{-t_{k+1}A_{g,t,{T\over t}}^2}\right]dt_1\cdots dt_k\right| \\
 \leq C_{14}  t^{k-n}    \left\| \psi
e^{-{ 1\over 2(k+1)}A_{g,t,{T\over t}}^2} \right\|
\end{multline}
for some constant $C_{14}>0$.

Now since ${\rm Supp}(\psi)\subset M\setminus \cup_{x\in B}
B_x(\delta)$, by \cite[Proposition 15.1]{BZ1}, one sees that there
exist $C_{15},\ C_{16}>0$ such that any $0<t\leq c_7$, $T\geq 1$,
\begin{align}\label{9.23}
\int_M{\rm Tr}\left[\psi(x)S_{{1\over \sqrt{  k+1 }}t,{1\over
\sqrt{ k+1 }}{T\over t}}(x,x)\psi(x)\right]d{\rm vol}_x\leq
C_{15}\exp\left(-{C_{16}T\over t^2}\right).
\end{align}

From (\ref{8.9}), (\ref{9.22}) and (\ref{9.23}), one sees
immediately that there exists $C_{17}>0$, $C_{18}>0$ such that for
any $1\leq k\leq n-1$, $0<t\leq c_7$ and $T\geq 1$, one has
\begin{multline}\label{9.24}
\left|\int_{\Delta_k}{\rm Tr}_s\left[N  e^{-t_1A_{g,t,{T\over
t}}^2}C_{t,{T\over t}} e^{-t_2A_{g,t,{T\over t}}^2}\cdots
C_{t,{T\over t}}e^{-t_{k+1}A_{g,t,{T\over t}}^2}\right]dt_1\cdots dt_k\right| \\
 \leq C_{17}
e^{-C_{18}T } .
\end{multline}

From (\ref{9.3}), (\ref{9.21}) and (\ref{9.24}), one gets
(\ref{9.1}).

The proof of Theorem \ref{t3.9} is completed. \ \ Q.E.D.

\section{Euler structure and   the Burghelea-Haller conjecture}\label{s10}
\setcounter{equation}{0}

In this section we recall several symmetric bilinear  torsions
introduced by Burghelea-Haller \cite{BH1, BH2} which are defined
by using the Euler structure introduced by Turaev \cite{T}. We
then apply our main result, Theorem \ref{t3.1},  to prove   a
conjecture due to Burghelea and Haller \cite[Conjecture 5.1]{BH2}.

Some applications on comparisons of various torsions are also
included.

\subsection{Euler and coEuler structures}

 Let $M$ be a closed oriented smooth manifold, with
$\dim M=n$. We   assume the vanishing of the Euler-Poincar\'{e}
characteristics of $M$, that is, $\chi(M)=0$. The set of Euler
structures with integral coefficients, Eul($M; \textbf{Z}$),
introduced by Turaev \cite{T}, is an affine version of $H_{1}(M;
\textbf{Z})$.

Let $X\in\Gamma(TM)$ be a    non-degenerate vector field on $M$
which means $X : M \longrightarrow{TM}$ is transversal  to the
zero section. Denote its set of zeros by ${\rm zero}({X})$. For
every $x\in{\rm zero}({X})$, there is a well-defined   Hopf index
${\rm IND}_{X}(x)\in \{\pm{1}\}$.

Any  Euler structure can be represented by a pair $(X, c)$ where
$c\in{C_{1}^{\rm sing}(M; \textbf{Z})}$ is a singular 1-chain
satisfying
\begin{align}\label{10.1}
\partial{c}={\rm e}(X):=\sum_{x\in{\rm zero}({X})}{\rm IND}_{X}(x)x.
\end{align}
Since  $\chi(M)=0$, the existence of $c$ is clear.

\begin{lemma}\label{t10.1}{\rm (\cite[Lemma 2.1]{BH2})}
Let $M$ be a closed smooth manifold with $\chi(M)=0$, let
$\mathfrak{e}\in{\rm Eul}(M; {\bf Z})$ be an Euler structure, and
let $x_0\in{M}$ be a base point. Suppose X is a non-degenerate
vector field on M with   ${\rm zero}({X})\neq{\emptyset}$. Then
there exists a collection of paths $\sigma_x$, $\sigma_x(0)=x_0$,
$\sigma_x(1)=x$, $x\in{\rm zero}({X})$, so that
\begin{align}\label{10.2}
\mathfrak{e}=\left[X, \sum_{x\in{\rm zero}({X})}{\rm
IND}_X(x)\sigma_x\right].
\end{align}
\end{lemma}

The set of coEuler structures ${\rm Eul}^*(M;\textbf{C})$ is an
affine version of $H^{n-1}(M; \textbf{C})$.

Let $g^{TM}$ as before be a Riemannian metric on $M$ with  the
associated  Levi-Civita connection denoted by $\nabla^{TM}$.

 Any
coEuler structure can be represented by $(g^{TM},\alpha)$ for some
$\alpha\in\Omega^{n-1}(M)$ such that
\begin{align}\label{10.3}
d\alpha=e\left(TM,\nabla^{TM}\right),
\end{align}
 where $e(TM,\nabla^{TM})$ is the Euler form defined in
(\ref{2.27}). Since $\chi(M)=0$, the existence of $\alpha$ is
clear.

If $[X,c]\in{\rm Eul}(M;\textbf{Z})$  and $[g^{TM},\alpha]\in{\rm
Eul}^*(M;\textbf{C})$, we call $[g^{TM},\alpha]$ is dual to
$[X,c]$ if for any closed one form $\omega\in \Omega^1(M)$ which
vanishes in a neighborhood of ${\rm zero}({X})$,
\begin{align}\label{10.4}
\int_{M
}\omega\wedge\left(X^*\psi\left(TM,\nabla^{TM}\right)-\alpha\right)=\int_c\omega,
\end{align}
 where
$\psi(TM,\nabla^{TM}) $ is  the Mathai-Quillen current (\cite{MQ})
associated with $g^{TM}$ defined in \cite[Definition 3.6]{BZ1}.

For any $[X,c]\in{\rm Eul}(M;\textbf{Z})$ and $g^{TM}$, the
existence of $\alpha$ is proved in \cite{BH1, BH2}.

\subsection{A proof of the Burghelea-Haller conjecture}\label{s10.1}

We make the same gemteric assumptions as in Section \ref{s3}. We
also assume $\chi(M)=0$ as in the previous subsection.

Recall that we have the Thom-Smale cochain complex
$(C^*(W^u,F),\partial)$ associated to a Morse function $f$ and a
Riemannian metric $g^{TM}$ verifying conditions in Section
\ref{s3.1}.

 Let
$x_0\in{M}$ be a fixed base point.

Let $\mathfrak{e}$ be an Euler structure.

For every critical point $x\in{B}$ of $f$ choose a path $\sigma_x$
with $\sigma_x(0)=x_0$ and $\sigma_x(1)=x$ so that $ [\nabla f,
\sum_{x\in{B}}(-1)^{{\rm ind}(x)}\sigma_x]$ is a representative of
$\mathfrak{e}$ (cf. Lemma \ref{t10.1}).

Let $b_{x_0}$ be a nondegenerate symmetric bilinear form  on the
fiber $F_{x_0}$ over $x_0$. For $x\in{B}$ define a nondegenerate
symmetric bilinear form $b_x$ on $F_x$ by parallel transport of
$b_{x_0}$ along $\sigma_x$ with respect to $\nabla^F$. The
collection of symmetric bilinear forms $\{b_x\}_{x\in{B}}$ defines
a nondegenerate symmetric bilinear form on the Thom-Smale cochain
complex $(C^*(W^u,F),\partial)$, which in turn defines an induced
symmetric bilinear form on ${\rm det}\, H^*(C^*(W^u,F),\partial)$.

Since $\chi(M)=0$, one sees easily that the above induced
symmetric bilinear form on ${\rm det}\, H^*(C^*(W^u,F),\partial)$
 does not depend on the choices of
$\{\sigma_x\}_{x\in{B}}$, $x_0$ and  $b_{x_0}$. It depends only on
  $F$, $\mathfrak{e}$ and $\nabla f$. We call it the Milnor-Turaev
  symmetric bilinear torsion and denote it by $\tau_{F, \mathfrak{e}}^{\nabla f}$.

On the other hand, let $b^F$ be a nondegenerate symmetric bilinear
form on the flat vector bundle $F$.

 For any $\alpha\in
\Omega^{n-1}(M)$ such that
  $d\alpha={  e}(TM,\nabla^{TM})$, following  Burghelea and
  Haller \cite{BH1, BH2}, one defines
\begin{align}\label{10.5}
\tau_{F, g^{TM}, b^F, \alpha}^{\rm an}=
 b_{(M,F,g^{TM},b^F)}^{RS} \cdot\exp\left(\int_{M}\theta\left({F,b^F}\right)\wedge\alpha\right)
\end{align}
and call it the Burghelea-Haller symmetric bilinear  torsion.

 By \cite[Theorem 4.2]{BH2}, we know that $\tau_{F, g^{TM}, b^F, \alpha}^{\rm
 an}$ does not depend on the choice of $g^{TM}$ and the smooth
 deformations of $b^F$. Thus we now denote it by $\tau_{F,   b^F, \alpha}^{\rm
 an}$.

We can now state the following equivalent version of the
Burghelea-Haller conjecture \cite[Conjecture 5.1]{BH2}.

\begin{thm}\label{t10.2}
  If $\mathfrak{e}= [\nabla f,
\sum_{x\in{B}}(-1)^{{\rm ind}(x)}\sigma_x]$ and  $(g^{TM},\alpha)$
are dual in the sense of (\ref{10.4}), then   we have
\begin{align}\label{10.6}
P^{\det H}_\infty\left(\tau_{F,   b^F, \alpha}^{\rm
an}\right)=\tau_{F, \mathfrak{e}}^{\nabla f}.
\end{align}
\end{thm}
{\it Proof}.  By \cite[Theorem 4.2]{BH2}, we may well assume that
$b^F$ is flat near $B$. Then $\theta(F,\nabla^F)=0$ near $B$.

By (\ref{10.4}), one has
\begin{align}\label{10.7}
\int_{M
}\theta\left(F,b^F\right)\left(X^*\psi\left(TM,\nabla^{TM}\right)
-\alpha\right)=\int_{c}\theta\left(F,b^F\right),
\end{align}
where $c=\sum_{x\in{B}}(-1)^{{\rm ind}(x)}\sigma_x$.

 From  Theorem 3.1 and (\ref{10.7}), we have in noting $X=\nabla
 f$ that,
\begin{multline}\label{10.8}
P_\infty^{{\rm det}H}\left(\tau_{F,  b^F, \alpha}^{\rm
an}\right)=P_\infty^{{\rm det}H}\left(b_{(M,F,g^{TM},b^F)}^{\rm
RS}\right)\cdot
\exp\left(\int_{M}\theta\left({F,b^F}\right)\wedge\alpha\right)\\
=b^{\cal M}_{( {M},F,b^F,-X)} \cdot\exp\left(-
\int_M\theta\left(F,b^F\right)X^*\psi\left(TM,
\nabla^{TM}\right)\right)\cdot\exp\left(\int_{M}\theta\left({F,b^F}\right)\wedge\alpha\right)\\
=b^{\cal M}_{(
{M},F,b^F,-X)}\cdot\exp\left(\int_{M}\theta\left({F,b^F}\right)
\wedge\left(\alpha-X^*\psi\left(TM,\nabla^{TM}\right)\right)\right)\\
 =b^{\cal M}_{(
{M},F,b^F,-X)}\cdot\exp\left(-\int_{c}\theta\left(F,b^F\right)\right).
\end{multline}
By \cite[(46)]{BH2}, we have
\begin{align}\label{10.9}
\tau_{F, \mathfrak{e}}^{X}=b^{\cal M}_{(
{M},F,b^F,-X)}\cdot\exp\left(-\int_{c}\theta\left(F,b^F\right)\right).
\end{align}

By (\ref{10.8}) and (\ref{10.9}), we get (\ref{10.6}).

The proof of Theorem \ref{t10.2} is completed. \ \ Q.E.D.

$\ $

\begin{cor}\label{t10.3} The Burghelea-Haller torsion $\tau_{F, b^F, \alpha}^{\rm
an}$ does not depend on   $b^F$ and the representative $\alpha$.
\end{cor}

\begin{Rem}\label{t10.4}  In view of the remarks given in
\cite[Section 7.3]{BH1}, Theorem \ref{t10.2} provides an
analytic interpretation of the Alexander polynomial in knot
theory.
\end{Rem}

\subsection{Comparison of $ b^{\rm RS}_{(M,F,g^{TM},b^{F})} $ with the usual
Ray-Singer torsion}

We still assume  $\chi(M)=0$.

Let $g^F$ be a Hermitian metric on $F$. Then one can construct the
Ray-Singer analytic torsion  as an inner product on $\det
H^*(M,F)$ (or equivalently as a metric on the determinant line,
cf. \cite{BZ1}). We denote the resulting inner product by $b^{\rm
RS}_{(M,F,g^{TM},g^{F})} $.

In this section, we prove the following comparison result between
$b^{\rm RS}_{(M,F,g^{TM},b^{F})} $ and $b^{\rm
RS}_{(M,F,g^{TM},g^{F})} $, which is also a consequence of
\cite[(5.13)]{BK3} and \cite[Theorem 1.4]{BK4}.

It is clear that the absolute value of the ratio of the symmetric
bilinear form and the inner product is well-defined.

\begin{prop}\label{t10.5} If $\dim M$ is odd, then the following identity holds,
\begin{align}\label{10.10}
\left|{  b_{(M,F,g^{TM},b^F)}^{\rm RS} \over
b_{(M,F,g^{TM},g^F)}^{\rm RS} }\right|=1.
\end{align}
\end{prop}

{\it Proof}. Let $\mathfrak{e}$ be an Euler class associated to
$\nabla f$ in the sense of Lemma \ref{t10.1}. Let $T^{\nabla
f}_{F, \mathfrak{e}}$ be the Redemeister inner product induced
from the Euler structure $\mathfrak{e}$. Then one verifies easily
that
\begin{align}\label{10.11}
\left|{\tau^{\nabla f}_{F, \mathfrak{e}}\over T^{\nabla f}_{F,
\mathfrak{e}}}\right|=1.
\end{align}

Let $[g^{TM},\alpha]$, $\alpha\in \Omega^{n-1}(M)$, be dual to the
Euler structure $\mathfrak{e}$ in the sense of (\ref{10.4}).

From (\ref{10.5}), (\ref{10.6}), (\ref{10.11}) and \cite[Theorem
0.2]{BZ1}, one deduces that
\begin{align}\label{10.12}
\left|{  b_{(M,F,g^{TM},b^F)}^{\rm RS} \over
b_{(M,F,g^{TM},g^F)}^{\rm RS}
}\right|=\left|\exp\left(\int_{M}\left(\theta\left(F,g^F\right)
-\theta\left(F,b^F\right)\right)\wedge\alpha\right)\right|.
\end{align}

Note that the left hand side of (\ref{10.12}) does not depend on
the Euler structure $\mathfrak{e}$.

By choosing different Euler structures, one sees that  for any
real closed form $\gamma \in \Omega^{n-1}(M)$ whose image    in
$H^*(M,{\bf R})$ lies in $H^*(M,{\bf Z})$, one has
\begin{align}\label{10.13}
{\rm Re}\left(\int_{M}\left(\theta\left(F,g^F\right)
-\theta\left(F,b^F\right)\right)\wedge\gamma\right)=0.
\end{align}
Then it is easy to see that (\ref{10.13}) should also hold for any
real closed form $\gamma\in\Omega^{n-1}(M)$. As a consequence, we
get the following equality in $H^1(M,{\bf R})$
\begin{align}\label{10.14}
{\rm Re}\left[ \theta\left(F,b^F\right)\right]=
\left[\theta\left(F,g^F\right)\right] .
\end{align}

Since $\dim M$ is odd implies $e(TM,\nabla^{TM})=0$, by
(\ref{10.3}), (\ref{10.12}) and (\ref{10.14}),  we get
(\ref{10.10}).

 The proof of Proposition
\ref{t10.5} is completed.\ \ Q.E.D.

\begin{Rem}\label{10.6}
In the  general case that $\dim M$ need not be odd, by the
consideration in the proof of \cite[Theorem 5.9]{BH2}, one sees
that there exists an anti-linear involution $J^F: F\rightarrow F$
such that
\begin{align}\label{10.15}
\left(J^F\right)^2={\rm Id}_F,\ \ \ b^F\left(
J^Fu,v\right)=\overline{b^F\left( u,J^Fv\right)},\ \ \ b^F\left(
u,J^Fu\right)\geq 0,\ \ \ u,\ v\in F.
\end{align}
Then
\begin{align}\label{10.16}
 g^F(u,v):= b^F\left(
u,J^Fv\right) ,\ \ \ u,\ v\in F,
\end{align}
defines a Hermitian metric on $F$. From (\ref{10.16}), we get
\begin{align}\label{10.17}
\left(g^F\right)^{-1}\nabla^Fg^F=
\left(J^F\right)^{-1}\left(\left(b^F\right)^{-1}\nabla^Fb^F\right)J^F+
\left(J^F\right)^{-1}\nabla^FJ^F.
\end{align}
From (\ref{10.17}), one gets
\begin{align}\label{10.18}
\theta\left(F,b^F\right)= \theta\left(F,g^F\right)- {\rm Tr}\left[
\left(J^F\right)^{-1}\nabla^FJ^F\right],
\end{align}
from which we get
\begin{align}\label{10.19}
{\rm Re}\left(\theta\left(F,b^F\right)\right)=
\theta\left(F,g^F\right) ,
\end{align}
which provides a direct proof of (\ref{10.14}).
\end{Rem}

\comment{Now with $b^F$, $g^F$, $J^F$ in the current situation,
since by (\ref{10.15}) one has
\begin{align}\label{10.20}
\det\left(J^F\right)=\pm 1 ,
\end{align}
one sees easily that the following identity holds,
\begin{align}\label{10.21}
 \left|{b^{\cal M}_{(
{M},F,b^F,-X)}\over b^{\cal M}_{( {M},F,g^F,-X)}}\right|=1,
\end{align}
where $b^{\cal M}_{( {M},F,g^F,-X)}$ is the Reidemeister inner
product on $\det H^*(C^*(W^u,F),\partial  )  $ induced by $g^F$.

From  (\ref{10.19}), (\ref{10.21}), Theorem \ref{t3.1} and
\cite[Theorem 0.2]{BZ1}, we get the following result which holds
with out the assumption that $\chi(M)=0$.

\begin{thm}\label{t10.7} For any nondegenerate symmetric bilinear form $b^F$ on $F$,
there exists a Hermitian metric $g^F$ on $F$ such that
(\ref{10.10}) holds.
\end{thm}}

\subsection{On Braverman-Kappeler's approach}

When $M$ is of odd dimension, Braverman-Kappeler
\cite{BK1}-\cite{BK4} developed another approach of complex valued
analytic torsions. In particular, a comparison result between the
analytic torsions defined in \cite{BK2} and \cite{BH2} is proved
in \cite[Theorem 1.4]{BK4}.

Here we point out that with the help of Theorem \ref{t10.2}, one
can identify, at least up to $\pm$, a locally constant defined in
\cite[(5.11)]{BK3} on the moduli space of representations of the
fundamental group $\pi_1(M)$. Indeed, this follows from
\cite[(5.4), (5.5)]{BK4} and the proved Burghelea-Haller
conjecture (cf. \cite[Conjecture 1.9 and Theorem 1.10]{BK4})
immediately. We leave the details to the interested readers.


\begin {thebibliography}{15}

\bibitem[ABP]{ABP} M. F. Atiyah, R. Bott and V. K. Patodi, On the
heat equation and the index theorem. {\it Invent. Math.} 19
(1973), 279-330.


\bibitem[BGS]{BGS} J.-M. Bismut, H. Gillet and C. Soul\'e,
Analytic torsions and holomorphic determinant line bundles I. {\it
Commun. Math. Phys.} 115 (1988), 49-78.



\bibitem[BL]{BL} J.-M. Bismut  and G. Lebeau, Complex immersions
and Quillen metrics. {\em Publ. Math. IHES.} 74 (1991), 1-297.

\bibitem[BZ1]{BZ1} J.-M. Bismut and W. Zhang, {\it An Extension of a
Theorem by Cheeger and M\"uller}. {\it Ast\'erisque} Tom. 205,
Paris, (1992).

 \bibitem[BZ2]{BZ2} J.-M. Bismut and W. Zhang,
Milnor and Ray-Singer metrics on the equivariant
determinant of a flat vector bundle. {\it Geom.  Funct. Anal.}
4 (1994), 136-212.

\bibitem[BrK1]{BK1} M. Braverman and T. Kappeler, A refinement of
the Ray-Singer torsion. {\it C. R. Acad. Sci. Paris}, 341 (2005),
497-502.

\bibitem[BrK2]{BK2} M. Braverman and T. Kappeler, Refined analytic
torsion as an element of the determinant line. {\it Preprint,}
math.DG/0510532.

\bibitem[BrK3]{BK3} M. Braverman and T. Kappeler, Ray-Singer type
theorem for the refined analytic torsion. {\it Preprint,}
math.DG/0603638.

\bibitem[BrK4]{BK4} M. Braverman and T. Kappeler, Comparison of
the refined analytic and the Burghelea-Haller torsions. {\it
Preprint,} math.DG/0606398.




 \bibitem [BuH1]{BH1} D. Burghelea and  S. Haller,  Torsion, as function on the space of representations.
  {\it Preprint,} math.DG/0507587.

 \bibitem [BuH2]{BH2} D. Burghelea and  S. Haller, Complex valued
 Ray-Singer torsion. {\it Preprint,} math.DG/0604484.





 \bibitem[C]{C} J. Cheeger, Analytic torsion and the heat equation.
 {\it Ann. of Math.} 109 (1979), 259-332.

 \bibitem[CH]{CH} S.   Chern and X.   Hu, Equivariant Chern
 character for the invariant Dirac operator. {\it Michigan Math.
 J.} 44 (1997), 451-473.

\bibitem[FT]{FT} M. Farber and V. Turaev, Poincar\'e-Reidemeister
metric, Euler structures and torsion. {\it J. Reine Angew. Math.}
520 (2000), 195-225.

\bibitem[Fe]{F} H. Feng, A remark on the noncommutative Chern
character. {\it Acta Math. Sinica} 46 (2003), 57-64. (in Chinese)

\bibitem[G1]{G1}
 E. Getzler, Pseudodifferential operators on supermanifolds and
 the Atiyah-Singer index theorem. {\it Commun. Math. Phys.} 92
 (1983), 163-178.

\bibitem[G2]{G2}
 E. Getzler, A short proof the local Atiyah-Singer index theorem.
 {\it Topology} 25 (1986), 111-117.


 \bibitem [HS]{HS} B. Helffer and J. Sj\"ostrand, Puis multiples
 en m\'ecanique semi-classique IV: Etude du complexe de Witten.
 {\it Comm. PDE} 10 (1985), 245-340.

\bibitem [KM]{KM} F. F. Knudson and D. Mumford, The projectivity
of the moduli space of stable curves I: Preliminaries on ``det''
and ``div''. {\it Math. Scand.} 39 (1976), 19-55.

\bibitem [L]{La} F. Laudenbach, On the Thom-Smale complex.
Appendix in [BZ1].





\bibitem[MQ]{MQ}   V. Mathai and D. Quillen, Superconnections,
Thom classes, and equivariant differential forms. {\it Topology}
25 (1986), 85-110.


 \bibitem[Mi]{Mi} J. Milnor, Whitehead torsion. {\it  Bull.
 Amer.  Math. Soc.} 72 (1966), 358-426.

 \bibitem[Mu1]{Mu1}  W. M\"uller,  Analytic torsion  and the  R-torsion
 of Riemannian manifolds. {\it Adv. in Math.} 28 (1978), 233-305.

  \bibitem[Mu2]{Mu2}  W. M\"uller,  Analytic torsion  and the  R-torsion
  for unimodular representations. {\it J. Amer. Math. Soc.} 6
  (1993), 721-753.


  \bibitem [Q1]{Q}  D. Quillen, Determinants of Cauchy-Riemann
  operators over a Riemann surface. {\it Funct. Anal. Appl.} 14
  (1985), 31-34.

  \bibitem [Q2]{Q1}  D. Quillen, Superconnections and the Chern
  character. {\it Topology} 24 (1985), 89-95.

\bibitem [RS]{RS} D. B. Ray and I. M. Singer, $R$-torsion and the
Laplacian on Riemannian manifolds. {\it Adv. in Math.} 7 (1971),
145-210.

\bibitem [S1]{S1} M. A. Shubin, Semiclassical asymptotics on
covering manifolds and Morse inequalities. {\it Geom. Funct.
Anals.} 6 (1996), 370-409.

\bibitem [S2]{S} M. A. Shubin, {\it Pseudodifferential Operators
and Spectral Operator}. Springer-Verlag, Berlin, 2001.

\bibitem [Si]{Si} B. Simon, {\it Trace Ideals and their
Applications}. London Mathematical Society, Lecture Notes Series,
35, Cambridge University Press, 1979.

\bibitem [Sm]{Sm} S. Smale, On gradient dynamical systems. {\it Ann.
of Math.} 74 (1961), 199-206.

\bibitem [SuZ]{SZ} G. Su and W. Zhang,  A Cheeger-M\"uller theorem
for complex valued Ray-Singer torsion. {\it Preprint,} 2006.

\bibitem[T]{T} V. Turaev, Euler structures, nonsingular vector
fields, and Reidemeister-type torsion. {\it Math. USSR-Izv.} 34
(1990), 627-662.

 \bibitem[W]{W}  E. Witten, Supersymmetry and Morse theory. {\it
 J. Diff. Geom.} 17 (1982), 661-692.


 \bibitem[Z]{Z} W. Zhang, {\it Lectures on Chern-Weil Theory and Witten Deformations,}
      Nankai Tracts in Mathematics, Vol. 4.         World Scientific, Singapore, 2001.




\end{thebibliography}
\end{document}